\documentclass[11pt]{article}
\pdfoutput=1
\newcommand{\inc}{.}

\usepackage[english]{babel}
\usepackage{amsmath,amssymb,amsfonts}
\usepackage[pdftex]{graphicx,color}
\usepackage{subfigure}
\usepackage{empheq} 

\newtheorem{thm}{Theorem}[section]
\newtheorem{lem}{Lemma}[section]
\newtheorem{rmk}{Remark}[section]

\newcommand{\brk}[1]{\left(#1\right)}          
\newcommand{\BRK}[1]{\left\{#1\right\}}        
\newcommand{\Abs}[1]{\left\vert #1 \right\vert}        

\newcommand{\Norm}[1]{\left\Vert #1 \right\Vert} 
\newcommand{\Normt}[1]{\left\Vert #1 (t,\cdot)\right\Vert} 
\newcommand{\Norms}[1]{\left\Vert #1 (s,\cdot)\right\Vert} 
\newcommand{\Normtp}[1]{\left\Vert \left(#1\right) (t,\cdot)\right\Vert} 
\newcommand{\pd}[2]{\dfrac{\partial#1}{\partial#2}}
\newcommand{\deriv}[2]{\dfrac{d#1}{d#2}}
\newcommand{\pdd}[2]{\dfrac{\partial^2#1}{\partial#2^2}}
\newcommand{\into}{\int_\Omega}

\newcounter{step}
\newcommand{\newstep}[2]{
\refstepcounter{step}
\subparagraph*{Step \arabic{step}: #1.}
\addcontentsline{toc}{subsection}{Step \arabic{step}: #1}
\label{#2}
}

\def\uinf{u_\infty}
\def\tauinf{\tau_\infty}
\def\finf{f_\infty}
\def\Ld{{L^2}}
\def\Li{{L^\infty}}

\def\nequiv{\not{\!\!\equiv}}
\def\tt{\widetilde}

\newcommand{\dps}[1]  {\displaystyle{#1} }
\def\proof{ {{\bf Proof.~~}}}
\def\endproof{\hfill$\diamondsuit$\par\medskip}
\numberwithin{equation}{section}
\newcommand{\acknowledgements}[1]{\paragraph*{Acknowledgements.}#1}

\let\pa=\partial
\let\f=\frac
\let\Om=\Omega
\def\R{\mathop{\mathbb R\kern 0pt}\nolimits}

\newcommand{\beno}{\begin{eqnarray*}}
\newcommand{\eeno}{\end{eqnarray*}}
\def\kgp{\kappa}
\def\mm{m}
\def\veps{\epsilon}
\def\vveps{\alpha}
\def\np{n}

\def\Ldloc{L^2_{loc}}

\author{David Benoit$^{(1,2)}$, Lingbing He$^{(3)}$, Claude Le Bris$^{(1,2)}$ and Tony Leli\`evre$^{(1,2)}$\\
\footnotesize{(1) CERMICS, Ecole des
  Ponts (ParisTech), 6 \& 8 Av. B. Pascal,
  77455 Marne-la-Vall\'ee, France.}\\
\footnotesize{(2) INRIA Rocquencourt, MICMAC team, B.P. 105, 78153 Le
  Chesnay Cedex, France.}\\
\footnotesize{(3) Department of Mathematical Sciences, Tsinghua University
Beijing 100084,  P. R.  China.} \\
\{benoitd,lebris,lelievre\}@cermics.enpc.fr ~~~ lbhe@math.tsinghua.edu.cn
}
\title{Mathematical analysis of a one-dimensional model for an aging 
 fluid}

\begin{document}
\maketitle

\begin{abstract}
  We study mathematically a system of partial differential equations
  arising in the modelling of an aging fluid, a particular class of non
  Newtonian fluids. We prove well-posedness of the equations in
  appropriate functional spaces and investigate the longtime behaviour
  of the solutions.
\end{abstract}

\section{Introduction}
Our purpose is to study mathematically a system of partial differential
equations arising in the modelling of some particular non Newtonian
fluids. These fluids are often called \emph{aging fluids}. Two physical
phenomena are indeed permanently competing within the flow of such
fluids. On the one hand, the fluid \emph{ages} in the sense that it
solidifies. On the other hand, aging is counterbalanced by a
flow-induced rejuvenation.

The specific modelling we consider has been proposed in \cite{derec-01}
on the basis of phenomenological arguments and experimental
observations. A coefficient $f$, called the \emph{fluidity} encodes
aging for all times and at every location within the fluid. The fluid is
solid where $f=0$, and behaves all the more as a liquid when $f$ grows.
Our mathematical study aims to contribute to better understand how well
such a model captures the essential phenomena at play in fluid aging.

For our study, we proceed in a one-dimensional setting corresponding
physically to the consideration of a laminar Couette flow. Our three
unknown fields, the velocity $u$, the shear stress $\tau$ and the
fluidity $f$ are defined as functions of a space variable $y$ varying in
the interval $[0,1]$. They are also, of course, functions of the time
$t\geq 0$. The specific system we choose for our study reads
\begin{subequations} \label{eq:couette}
  \begin{empheq}[left=\empheqlbrace]{align}
    \rho \pd{u}{t} &= \eta \pdd{u}{y}+ \pd{\tau}{y},\label{eq:couettea}\\
    \lambda \pd{\tau}{t} &= G \pd{u}{y} - f \tau, \label{eq:couetteb}\\
    \pd{f}{t} &= (-1+\xi |\tau|) f^2-\nu f^3. \label{eq:couettec}
  \end{empheq}
\end{subequations}

Six dimensionless coefficients, all positive, constant in time and
throughout the domain, are present in the system: the density $\rho$,
the viscosity $\eta$, the characteristic relaxation time~$\lambda$, the
elastic modulus $G$, and two coefficients $\xi$ and $\nu$ specifically
related to the equation for the evolution of the fluidity $f$.  System
\eqref{eq:couette} is a fully coupled system of three equations. The
first two equations are classical in nature. The first one is the
equation of conservation of momentum for $u$. The second equation rules
the evolution of the shear stress $\tau$. The non classical ingredient
therein (as opposed, say, to an Oldroyd-B type equation) is the presence
of an extra parameter, the fluidity $f$, the role of which is formally
similar to that of an inverse time in a relaxation phenomenon. We note
that when $f$ is a constant in time and throughout the domain, the
equation agrees with the one-dimensional Oldroyd-B equation considered
{\em e.g.} in~\cite{Boyaval2009,GUILLOPE1990}. The third equation is of
the form of one of the many such evolution equations suggested
in~\cite{derec-01}. We hope it is, in this respect, prototypical of a
general class. It models the evolution of the fluidity $f$ in function
of the stress tensor.  The right-hand side of~\eqref{eq:couettec} may
differ from one model to another. The important ingredient is the
presence of two competitive terms: a negative term modelling aging and a
positive term modelling rejuvenation. For mathematical convenience, we
have taken two particular instances of these two terms.

We examine well-posedness and longtime behaviour for
system~\eqref{eq:couette}. Because we provide self-contained proofs, our
study is rather long. Similar questions on a different, although related
model for a viscoelastic fluid, have been examined
in~\cite{Renardy2009e,Renardy2009}.  \newline

Our article is articulated as follows.

To start with, we prove in Section~2 that the system under consideration
admits a global-in-time solution in appropriate functional spaces. The
solution is shown to be unique, and indeed strong. System
\eqref{eq:couette} is thus satisfied in a classical sense. Our precise
statement is the object of Theorem~\ref{existence}. The bulk of
Section~2 consists of our proof. The arguments are standard arguments of
mathematical fluid dynamics: formal \emph{a priori} estimates,
approximation, rigorous \emph{a priori} estimates, convergence.  The
many nonlinearities present in system~\eqref{eq:couette} however prevent
us, in the current state of our understanding, from extending our
analysis to settings in dimensions higher than or equal
to~two. Technically, this is related to the fact we repeatedly use, in
our arguments, that $H^1$ functions are $L^\infty$ functions, a
specificity of the one-dimensional setting of course.

In Sections 3 and 4, we study the long time behaviour of the
solution. Section~3 deals with return to equilibrium. We supply the
system with \emph{homogeneous} Dirichlet boundary conditions for the
velocity and investigate whether the flow converges to a steady state.
For homogeneous Dirichlet boundary conditions, the steady states are
$(u\equiv 0, \tau\equiv c, f\equiv 0),$ where~$c$ is a constant
throughout the domain.  The long-time convergence to these steady-states
sensitively depends, in system~\eqref{eq:couette}, of the fluidity $f$.
The situation is qualitatively different depending on the fluidity $f_0$
at initial time.  The more delicate, but of course more interesting,
case mathematically is the case where the fluidity $f_0$ at initial time
does not vanish everywhere: a part of the material, possibly the whole
of it, is originally fluid.  Section~\ref{sec:longhom1} addresses this
case.  We show (and the proof is quite substantial even in the
one-dimensional setting we consider) that the flow converges to the null
steady-state in suitable functional norms. The precise statement is the
purpose of Theorem~\ref{th:longhom1}. The convergence is then shown to
be polynomial in time, for all three fields $u$, $\tau$ and $f$. The
rates of convergence are made precise in
Theorem~\ref{th:longhomt}. Numerical simulations we perform in Section~5
will show these rates are indeed sharp.  It is interesting to emphasize
the physical signification of our mathematical results. With regard to
modelling, the convergence of the fluidity~$f$ to zero that we
establish, under homogeneous Dirichlet boundary conditions, means that
when left at rest, the fluid progressively solidifies, a certainly
intuitive fact. In addition, for $u$ and $\tau$, the rate of convergence
sensitively depends on the size of the region where, originally, the
material is liquid (a size measured by our parameter~$\beta$ defined in
\eqref{def:beta} and present in the right-hand sides of the estimates of
Theorem~\ref{th:longhomt}). The larger the liquid region the quicker the
convergence of both the velocity and the shear stress to zero. It is not
completely clear to us whether the latter qualitative behaviour is or
not compatible with experimental observations or physical intuition.

If the material is entirely solid at initial time, that is $f_0\equiv 0$
everywhere, the behaviour is quite different. Then the material stays
solid for all times, while the velocity and shear stress vanish
exponentially fast. We present the simple analysis of this behaviour in
Section~\ref{sec:longhom0}. Note that the result agrees with simple
physical intuition.

Non-homogeneous boundary conditions, studied throughout
Section~\ref{sec:nhom}, are, as always for questions related to
long-time behaviours, significantly more intricate to address. We adopt
constant boundary conditions, respectively $u=0$ and $u=a>0$ at $y=0$
and $y=1$. We begin by showing in Section~\ref{ssec:ss} that, when we
impose that the fluidity is strictly positive everywhere, there exists a
\emph{unique} steady state.  We next show in Section~\ref{ssec:stab}
that this steady state is stable under small perturbations. our precise
result is stated in Theorem~\ref{th:stab}. When the perturbations of the
state are not small, analyzing return to equilibrium is, in general,
beyond our reach. We are however able to show that, when we assume a
particular form of the initial condition (namely linear velocity,
constant shear stress, constant positive fluidity), then return to
equilibrium does hold true even if the initial condition is not close to
the steady state.  Some suitable assumptions relating the size of the
parameters in system~\eqref{eq:couette} and the non-zero boundary
condition $a$ are also needed (see condition~\eqref{as}). Our precise
result is Theorem~\ref{th:ode}. The reason why we have to assume this
specific form of the initial conditions is purely technical (and our
numerical simulations will actually show that these restrictions are, in
practice, unnecessary). In that case, system~\eqref{eq:couette} reduces
to a two-dimensional system of ordinary differential equations, for
which Poincar\'e-Bendixson Theory allows us to understand the longtime
behaviour. Our study is performed in Section~\ref{ssec:edo}.

As briefly mentioned above, Section~5 presents some numerical
simulations. We first show that the rates of convergence estimated by
our various mathematical arguments in the various regimes considered in
Sections~\ref{sec:hom} and~\ref{sec:nhom} are indeed sharp. We also
investigate numerically the stability of the steady state.  Our
simulations show that, irrespective of the size of the initial
perturbation (and thus in a more general regime than that for our
mathematical arguments), the fluid returns to equilibrium, or more
generally converges to the suitable steady state. The rates of
convergence are also examined.  \newline

We conclude this introduction by mentioning that, despite their
limitations, our results show that the model derived in~\cite{derec-01}
does adequately account for aging and rejuvenation. However, two
shortcomings need to be emphasized. Both originate from the mathematical
nature of equation~\eqref{eq:couettec} (and are actually related to the
fact that the Cauchy Lipschitz theory applies to this equation). First,
when $f$ vanishes, then $f$ remains zero for all subsequent times. This
property, present everywhere in our mathematical study, prevents
fluidification to occur after solidification. This clearly limits the
range of materials covered by the modelling (compare muds and concrete,
say).  Second, $f$ can only vanish asymptotically and never in finite
time unless it is already zero before. Otherwise stated, solidification
can occur, but never in finite time: again a modelling limitation. The
one usefulness, if any, of our study, is therefore to point out that a
mathematically well founded model where fluidification \emph{and}
solidification compete on an equal footing is still to be derived.  Our
study implicitly points to suitable directions to this end.

Further mathematical investigations on models for aging fluids will be
presented in~\cite{these-david}.

\section{Global existence and uniqueness}
In this section, we establish the following global existence and
uniqueness result for system~\eqref{eq:couette} supplied with initial
conditions~$u_0,\tau_0,f_0$ and the boundary conditions~$u(t,0)=0$
and~$u(t,1)=a\ge0$ for all time $t\in[0,T]$ (where~$a$ is a constant
scalar).

\begin{thm}\label{existence}
  Recall that~$\Omega$ is the \emph{one-dimensional} domain~$[0,1]$ and
  that~$T>0$ is fixed. Consider the initial data
  \begin{align}\label{as:f0nn}
    \brk{u_0,\tau_0,f_0} \in H^1(\Omega)^3\mbox{ with } f_0\ge0.
  \end{align}
  Then there exists a unique global solution~$(u,\tau,f)$ to
  system~\eqref{eq:couette} such that for any~$T>0$,
  \begin{align}\label{reg}
    \brk{u,\tau,f}& \in \brk{C([0,T];H^1)\cap L^2([0,T];H^2)}\times
    C([0,T];H^1) \times C([0,T];H^1)
  \end{align}
  and~$f\ge0$ for all~$x\in\Omega$ and~$t\in[0,T]$.\\
  In addition, we have
  \begin{align}\label{regdt}
    \brk{\pd{u}{t},\pd{\tau}{t},\pd{f}{t}}\in L^2([0,T];L^2) \times
    C([0,T];L^2)\times C([0,T];L^2),
  \end{align}
  so that the equations in~\eqref{eq:couette} are all satisfied in the
  strong sense in time.
\end{thm}

Before we get to the proof, we eliminate the non-homogeneous Dirichlet
boundary condition, introducing the auxiliary velocity field
$$u(t,y) - a y.$$
This velocity field, which we still denote by~$u$, solves the system
\begin{subequations} \label{tcoue}
  \begin{empheq}[left=\empheqlbrace]{align}
    \rho \pd{u}{t} &= \eta \pdd{u}{y} + \pd{\tau}{y} \label{tcouea},\\
    \lambda \pd{\tau}{t} &= G \pd{u}{y} - f \tau+Ga\label{tcoueb}, \\
    \pd{f}{t} &= (-1+\xi |\tau|) f^2-\nu f^3\label{tcouec},
  \end{empheq}
\end{subequations}
supplied with the \emph{homogeneous} Dirichlet boundary
conditions~$u(t,0)=0$ and~$u(t,1)=0$ for all time~$t\in[0,T]$ and
initial conditions~$u_0,\tau_0,f_0 \in H^1(\Omega)$. The proof of
Theorem~\ref{existence} will actually be completed on
system~\eqref{tcoue}.  The result on~\eqref{eq:couette} then immediately
follows.

\medskip\proof The proof falls in eight steps. The first five steps
consist in deriving \emph{formal} a priori estimates.  These estimates
are next made rigorous for a sequence of approximate solution in
Step~\ref{s:cons}. The convergence of this sequence is proven in
Step~\ref{s:conv}, thereby establishing existence of a solution
to~\eqref{tcoue}. Step~\ref{s:uniq} addresses uniqueness.
\newstep{Non-negativity of the fluidity}{s:pos} Let us first formally
prove that $f\geq0$. Fix~$y\in\Omega$ and introduce
\begin{align*}
  E_0=\BRK{y\in\Omega, f_0(y)>0}.
\end{align*}
For~$y\in \Omega \setminus E_0$, we have~$f_0(y)=0$ and thus~$f(t,y)=0$
for all time~$t\in[0,T]$ because of~\eqref{tcouec}. On the other hand,
for~$y\in E_0$, we now show that~$f(t,y)>0$ for all time~$t\in[0,T]$. We
argue by contradiction and suppose, by continuity of~$f(\cdot,y)$, that
\begin{align*}
  t_m=\inf\BRK{t\in(0,T],f(t,y)=0}< T.
\end{align*}
The Cauchy-Lipschitz Theorem applied to~\eqref{tcouec} with zero as
initial condition at time~$t_m$ implies that~$f(t,y)=0$
for~$t\in\brk{t_m-\varepsilon,t_m+\varepsilon}$ for~$\varepsilon>0$,
which contradicts the definition of~$t_m$.

We have therefore shown that~$f$ stays zero where it is zero, and stays
positive where it is positive, which in particular implies
non-negativity everywhere.

\newstep{Formal first energy estimates}{s:fee} We again argue
formally. We first multiply the evolution equation~\eqref{tcouea} on~$u$
by~$u$ itself and integrate over the domain. This gives a first estimate
\begin{align} \label{enu} \f12\rho\f{d}{d t}\Normt{u}_\Ld^2+
  \eta\Normt{\pd{u}{y}}_\Ld^2=\into \brk{\pd{\tau}{y}
    u}(t,\cdot).  \end{align} Similarly, we multiply the evolution
equation~\eqref{tcoueb} by~$\tau$ and integrate over~$\Omega$ to find
\begin{align} \label{entau} \f12\lambda\f{d}{d
    t}\Normt{\tau}_\Ld^2+\Normt{\sqrt{f}\tau}_\Ld^2=G\into
  \brk{\pd{u}{y} \tau} (t,\cdot) +G\,a\,\bar{\tau}(t),\end{align} where
we denote by
\begin{align} \label{eq:avg} \bar{q}(t) =\dps \into q(t,y) d y
\end{align}
the average over~$\Omega$ of a function~$q:(t,y)\in [0,T] \times \Omega \rightarrow \R$. \\
Combining estimates~\eqref{enu} and~\eqref{entau} and using integration
by parts and the fact that~$u$ vanishes on the boundary, we obtain
\begin{align}\label{en1} \f12\f{d}{d t}\brk{G\rho
    \Normt{u}_\Ld^2+\lambda\Normt{\tau}_\Ld^2}+G\eta\Normt{\pd{u}{y}}_\Ld^2+
  \Normt{\brk{\sqrt{f}\tau}}_\Ld^2=G\,a\,\bar{\tau}(t).\end{align} We
now turn to~\eqref{tcouec}. Integrating~\eqref{tcouec} over~$\Omega$
yields
\begin{align} \label{ineq:fl1} \f{d}{d
    t}\Normt{f}_{L^1}+\Normt{f}_\Ld^2+\nu\Normt{f}_{L^3}^3 = \xi \into
  \brk{|\tau| f^2} (t,\cdot).
\end{align}
The Young inequality \beno \xi |\tau| f^2 = \sqrt\nu f^{\frac32} \cdot
\frac\xi{\sqrt\nu} |\tau|f^{\f12}\le \frac\nu2 f^3 + \frac{\xi^2}{2\nu}
f\tau^2 \eeno then yields
\begin{align}\label{en2} \f{d}{d
    t}\Normt{f}_{L^1}+\Normt{f}_\Ld^2+\f{\nu}{2}\Normt{f}_{L^3}^3\le
  \frac{\xi^2}{2\nu} \Normt{\brk{\sqrt{f}\tau}}_\Ld^2. \end{align}
Collecting ~\eqref{en1} and~\eqref{en2}, we obtain
\begin{align} \f12\f{d}{d t}&\brk{G\rho
    \Normt{u}_\Ld^2+\lambda\Normt{\tau}_\Ld^2+\frac{2\nu}{\xi^2}
    \Normt{f}_{L^1} } \nonumber\\
  &+G\eta\Normt{\pd{u}{y}}_\Ld^2+ \f12 \Normt{\brk{\sqrt{f}\tau}}_\Ld^2
  \le C\,a\Normt{\tau}_\Ld^2 \label{ebb}
\end{align}
where~$C$, here and throughout our text, denotes a constant, the actual
value of which is independent from~$T$ and only depends on the
domain~$\Omega$ and
the coefficients~$\rho, \eta, \lambda, G, \xi, \nu$ in~\eqref{tcoue} .\\
Applying the Gronwall Lemma to~\eqref{ebb}, we obtain
\begin{align}\label{eb} \sup_{t\in [0,T]}
  \brk{\Normt{u}_\Ld^2+\Normt{\tau}_\Ld^2+\Normt{f}_{L^1}} +\int_0^T
  \brk{\Normt{u}_{H^1}^2 +\right.&\left.
    \Normt{\brk{\sqrt{f}\tau}}_\Ld^2 } dt\nonumber\\
  &\le C_{0,T},
\end{align}
where~$C_{0,T}$ is a constant depending not only on~$\Omega, \rho, \eta,
\lambda, G, \xi, \nu$, but also on the boundary condition~$a$, the
initial data~$u_0,\tau_0,f_0$ and the time~$T$. \\
\begin{rmk}\label{rmk:chom1}
  For \emph{homogeneous} boundary conditions, that is~$a=0$, we mention
  that the right-hand sides of~\eqref{en1} and thus~\eqref{ebb} vanish.
  The constant~$C_{0,T}$ in~\eqref{eb} therefore does not depend on~$T$
  and we get a bound uniform in time.
\end{rmk}

\newstep{A priori estimates on an auxiliary function}{s:apf} Denote
\beno g(t,y)=\int_0^y (\tau(t,x)-\bar{\tau}(t))dx.\eeno This
function~$g$ satisfies \emph{homogeneous} Dirichlet boundary conditions
and formally solves \beno \pdd{g}{y}=\pd{\tau}{y}. \eeno
Using~\eqref{tcouea} and~\eqref{tcoueb}, which respectively imply \beno
\rho \pd{u}{t}=\eta \frac{\partial^2}{\partial y^2} \brk{ u+\f1{\eta}g}
,\eeno and \beno \lambda\pd{g}{t}=-\int_0^y(f\tau-\overline{f\tau})d x+G
u,\eeno we remark that the auxiliary function
\begin{align} \label{def:U}
  U &= u + \f1\eta \int_0^y \brk{\tau-\bar{\tau}}\\
  &= u+\f1\eta g \nonumber
\end{align}
solves:
\begin{align} \label{eq:U} \pd{U}{t}= \f{\eta}{\rho}
  \pdd{U}{y}-\f{1}{\lambda\eta}\int_0^y(f\tau-\overline{f\tau})d
  x+\frac{G}{\lambda\eta}u. \end{align} Multiplying
equation~\eqref{eq:U} by~$\pdd{U}{y}$ and integrating over $\Omega$
yields
\begin{align*}
  \f12 \f{d}{dt}&
  \Normt{\pd{U}{y}}_\Ld^2+\f{\eta}{\rho}\Normt{\pdd{U}{y}}_\Ld^2
  \nonumber\\&\leq C\brk{ \Normt{\brk{f\tau}}_{L^1} \into
    \Abs{\pdd{U}{y}}(t,\cdot) +\into \Abs{u \pdd{U}{y}}(t,\cdot)},
\end{align*}
after elementary manipulations in the right-hand side.  Then using the
Young and the Cauchy-Schwartz inequalities, we obtain
\begin{align} \label{en3} \f12 \f{d}{d
    t}\Normt{\pd{U}{y}}_\Ld^2+\f{\eta}{2\rho}\Normt{\pdd{U}{y}}_\Ld^2&\leq
  C
  \brk{\Normt{f}_{L^1}\Normt{\brk{\sqrt{f}\tau}}_\Ld^2+\Normt{u}_\Ld^2}.
\end{align}
Using~\eqref{eb}, we know that the right-hand side is~$L^1(0,T)$. In
view of our regularity assumptions on the initial conditions, we
have~$\left. \pd U y \right\vert_{t=0}
=\pd{u_0}{y}+\dfrac1{\eta}(\tau_0-\overline{\tau_0}) \in L^2(\Omega)$.
We therefore deduce from~\eqref{en3} that
\begin{align} \label{reg:U} U\in L^\infty([0,T],H^1_0)\cap
  L^2([0,T],H^2).
\end{align}

\newstep{$L^\infty$ estimates}{s:li} We are now in position to obtain
(again formal) $L^\infty$-bounds on~$\tau$ and~$f$. We consider the
evolution equation~\eqref{tcoueb}, which we rewrite in terms of~$U$
defined by~\eqref{def:U} and using~$\bar\tau$ defined for~$\tau$ as
in~\eqref{eq:avg}: \beno \lambda \pd{\tau}{t} =
G\pd{U}{y}-\brk{f+\f{G}{\eta}}\tau+\f{G}{\eta}\bar{\tau}+Ga. \eeno
Multiplying this equation by~$\tau$, we obtain \beno \f{\lambda}2\f{d}{d
  t}|\tau|^2+\brk{f+\f{G}{\eta}}|\tau|^2\leq C\brk{ |\tau|
  \cdot\Abs{\pd{U}{y}} + |\tau|\cdot \Norm{\tau}_\Ld+a\,|\tau|},\eeno so
that, repeatedly applying the Young inequality,
\begin{align}\label{en4} \f{\lambda}2\f{d}{d
    t}|\tau|^2+\brk{f+\f{G}{2\eta}}|\tau|^2\leq
  C\brk{\Abs{\pd{U}{y}}^2+\Norm{\tau}_\Ld^2+a}.\end{align} We apply the
Gronwall Lemma to~\eqref{en4} and use~$\pd Uy\in L^2([0,T],\Li)$ because
of~\eqref{reg:U}, estimate~\eqref{eb} and~$\tau_0\in H^1(\Omega)$ to
obtain
\begin{align}\label{eq:linftau}
  \Normt{\tau}_\Li\leq C_{0,T}
\end{align}
that is,~$\tau\in L^\infty([0,T],L^\infty)$.
 
As for the function~$f$, using the Duhamel formula for the evolution
equation~\eqref{tcouec} rewritten as \beno \pd{f}{t} = (-f-\nu f^2) f +
\xi |\tau| f^2, \eeno we obtain, for almost all~$y\in\Omega$,
\begin{align*}
  f(t,y)&=e^{-\int_0^t (f+\nu f^2)(s,y)d s} f_0(y)+ \xi\int_0^t
  e^{-\int_s^t (f+\nu f^2)(t^\prime,y) d t^\prime}
  \brk{|\tau|f^2}(s,y)ds\\
  &\le f_0(y)+\f\xi\nu \Norm{\tau}_{L^\infty_T(L^\infty)} \int_0^t
  e^{-\int_s^t \nu f^2(t^\prime,y) \nu d t^\prime} \nu f^2(s,y)ds,
\end{align*}
where we have used the non-negativity of~$f$ and the previously
derived~$L^\infty$-bound on~$\tau$ to obtain the second line. The above
equation leads to
\begin{align}\label{en5}
  f(t,y)&\le f_0(y)+\f{\xi}{\nu}\Norm{\tau}_{L^\infty_T(L^\infty)}
  \brk{1-e^{-\int_0^t \nu f^2(s,y) d s}}
  \nonumber\\
  &\le
  f_0(y)+\f{\xi}{\nu}\Norm{\tau}_{L^\infty_T(L^\infty)}. \end{align}
Using that~$f_0\in H^1$ and that we work in a one-dimensional setting,
we obtain that~$f\in L^\infty([0,T],L^\infty)$.

\begin{rmk} \label{rmk:chom2} For \emph{homogeneous} boundary
  conditions, the Gronwall Lemma applied to~\eqref{en4} implies
  \begin{align} \label{eq:tauinfh} \Normt{\tau}_\Li^2 \le
    \Norm{\tau_0}_\Ld^2 e^{-\f{G}{\lambda\eta}t} + \int_0^t \Norm{\pd
      Uy(s,\cdot)}_\Li^2 ds + \sup_{t\in[0,T]} \Normt{\tau}_\Ld^2
    \int_0^t e^{\f{G}{\lambda\eta}(s-t)} ds.
  \end{align}
  Moreover, as explained at the end of Step~\ref{s:fee}, the
  constant~$C_{0,T}$ in~\eqref{eb} does not depend on~$T$. Hence, the
  right-hand side of~\eqref{en3} and the bound
  in~$L^2([0,T];L^\infty)$-norm for~$\pd Uy$, deduced
  from~\eqref{reg:U}, also do not depend on~$T$. It follows
  from~\eqref{eq:tauinfh} that the~$L^\infty$-bound~\eqref{eq:linftau}
  on~$\tau$ is uniform in time. Equation~\eqref{en5} yields a similar
  conclusion for the bound on~$f$.
\end{rmk}

\newstep{Second a priori estimates}{s:see} In order to get estimates on
higher order derivatives, we now differentiate with respect to~$y$ the
evolution equation~\eqref{tcoueb} and obtain
\begin{align} \lambda \frac{\partial}{\partial t}\brk{\pd{\tau}{y}}&= G
  \pdd{u}{y}-f\pd{\tau}{y}-\pd{f}{y} \tau \nonumber\\ &= G
  \pdd{U}{y}-\f{G}{\eta}\pd{\tau}{y}-f\pd{\tau}{y}-\pd{f}{y}\tau. \label{tcouebd}\end{align}
Likewise, we differentiate with respect to~$y$ the evolution
equation~\eqref{tcouec} and get
\begin{align} \label{tcouecd} \frac{\partial}{\partial
    t}\brk{\pd{f}{y}}= \xi \pd{|\tau|}{y}
  f^2+2(\xi|\tau|-1)f\pd{f}{y}-3\nu f^2\pd{f}{y}. \end{align}
Multiplying equations~\eqref{tcouebd} and~\eqref{tcouecd} respectively
by~$\pd{\tau}{y}$ and~$\pd{f}{y}$, integrating over the domain, summing
up and using that both~$\tau$ and~$f$ are in~$L^\infty([0,T],L^\infty)$,
we obtain
\begin{align*} \f12\f{d}{d
    t}&\brk{\lambda\Normt{\pd{\tau}{y}}_\Ld^2+\Normt{\pd{f}{y}}_\Ld^2}
  \\& \leq C_{0,T} \into \brk{ \pdd{U}{y}\pd{\tau}{y} +
    \brk{\pd{\tau}{y}}^2 + \pd{f}{y} \pd{\tau}{y} + \pd{|\tau|}{y}
    \pd{f}{y} + \brk{\pd{f}{y}}^2 }(t,\cdot).
\end{align*}
Repeatedly applying the Young inequality and using that $\tau$ belongs
to~$L^2([0,T],H^1)$ (which implies~$\Abs{\pd{|\tau|}{y}(t,y)}=
\Abs{\pd{\tau}{y}(t,y)}$ for almost all~$t\in[0,T]$ and~$y\in\Omega$),
we obtain
\begin{align}
  \f12\f{d}{d
    t}&\brk{\lambda\Normt{\pd{\tau}{y}}_\Ld^2+\Normt{\pd{f}{y}}_\Ld^2}
  \nonumber\\& \leq C_{0,T}
  \brk{\Normt{\pd{\tau}{y}}_\Ld^2+\Normt{\pd{f}{y}}_\Ld^2+\Normt{\pdd{U}{y}}_\Ld^2
  }. \label{eq:sape}\end{align} We apply the Gronwall Lemma
to~\eqref{eq:sape}, use that~$\tau_0,f_0\in H^1(\Omega)$ and the
estimate~\eqref{reg:U} to obtain that~$\tau, f\in L^\infty([0,T],H^1)$.

It follows from the definition~\eqref{def:U} and the
estimate~\eqref{reg:U} that~$u\in L^\infty([0,T],H^1)\cap
L^2([0,T],H^2)$.

\newstep{Construction of an approximate solution}{s:cons} Now that we
have established all the necessary formal a priori estimates, we turn to
the construction of a sequence of approximating solutions
to~\eqref{tcoue} on which we will rigorously derive these a priori
estimates. We introduce, for~$n\ge1$, the sequence of systems
\begin{subequations}\label{eq:atcoue}
  \begin{empheq}[left=\empheqlbrace]{align}
    & \rho \pd{u_n}{t} = \eta \pdd{u_n}{y} + \pd{\tau_n}{y} \label{eq:atcouea},\\
    & \lambda \pd{\tau_n}{t} = G \pd{u_n}{y} - f_{n-1} \tau_n+Ga \label{eq:atcoueb},\\
    & \pd{f_n}{t} = (-1+\xi |\tau_n|) f_{n-1} f_n-\nu f_{n-1}
    f_n^2 \label{eq:atcouec},
  \end{empheq}
\end{subequations}
supplied with the \emph{homogeneous} Dirichlet boundary conditions
$u_n(t,0)=0$ and~$u_n(t,1)=0$ for all time~$t\in[0,T]$ and initial
conditions~$(u_{n0},\tau_{n0},f_{n0})=(u_0,\tau_0,f_0)$. We actually use
the initial condition $(u_0,\tau_0,f_0)$ also to initialize the
iterations in~$n$, thus the co\"incidence of notation.

We argue by induction. Consider
\begin{align*}
  \brk{u_{n-1},\tau_{n-1},f_{n-1}}& \in \brk{C([0,T];H^1)\cap
    L^2([0,T];H^2)}\times C([0,T];H^1) \times C([0,T];H^1)
\end{align*}
and~$f_{n-1}\ge0$.  We first show that there exists a unique
solution~$(u_n, \tau_n, f_n)$ to~\eqref{eq:atcoue} belonging to the same
functional spaces and such that~$f_n \ge 0$.  For this purpose, we
decompose~\eqref{eq:atcoue} into two subsystems: the linear (Oldroyd-B)
type model coupling the evolution equations~\eqref{eq:atcouea} on~$u_n$
and~\eqref{eq:atcoueb} on $\tau_n$ on the one hand and the ordinary
differential equation~\eqref{eq:atcouec} on~$f_n$ satisfied for
all~$y\in\Omega$ on the other hand.  The existence and uniqueness of a
solution~$(u_n,\tau_n)$ in the space $\brk{C([0,T];H^1)\cap
  L^2([0,T];H^2)}\times C([0,T];H^1)$ for the former system is obtained
using a classical approach (see for instance~\cite{GUILLOPE1990} for a
very close system).
We now turn to~$f_n$. We show that~$f_n$ exists in~$C([0,T];H^1)$
and~$f_n\ge0$. The equation~\eqref{eq:atcouec} writes
\begin{align} \label{eq:fngen}
  \pd{f_n}t&=\psi(t,f_n,y),\\
  \left. f_n \right\vert_{t=0}&=f_0\nonumber
\end{align}
where~$\psi$ is a function from~$[0,T]\times\R\times\Omega$ to~$\R$.

We first fix~$y\in\Omega$ and show that the function~$f_n(\cdot,y)$ is
continuous in time and non-negative.  The function~$\psi$ is continuous
in its first two variables and locally Lipschitz in its second
variable. The Cauchy-Lipschitz Theorem shows there exists a unique local
solution with~$f_0(y)$ as initial condition. Let~$[0,T^*)$ be the
interval of existence of the maximal solution for positive time.  For
all~$t\in[0,T^*)$, we have~$f_n\ge 0$, using Step~\ref{s:pos}.  In
addition, since~$f_{n-1}$ and~$f_n$ are both
non-negative,~\eqref{eq:atcouec} implies for all~$t\in[0,T^*)$,
\begin{align}
  \pd{f_n}t&\le \xi \Abs{\tau_n} f_{n-1} f_n \nonumber\\
  &\le \xi \Norm{\tau_n}_{C_T(\Li)} \Norm{f_{n-1}}_{C_T(\Li)}
  f_n, \label{eq:bofn}
\end{align}
using that both~$\tau_n$ and~$f_{n-1}$ belong to~$C([0,T];H^1)$. The
Gronwall Lemma then proves that~$f_n$ remains bounded on~$[0,T^*]$ and
thus we have established existence and uniqueness on~$[0,T]$.\\
We now turn to the local property of continuity of~$f_n$ as a function
of~$y$. We use that the function~$\psi$ is continuous in~$y$, because
both~$\tau_n$ and~$f_{n-1}$ are continuous in~$y$ in our one-dimensional
setting and the theorem on the continuous dependence on a parameter for
ordinary differential equations of the form~\eqref{eq:fngen} (see {\em
  e.g.}~\cite[Theorem~1.11.1, p.~126]{Cartan1967}).  We now show
that~$\pd{f_n}{y}$ belongs to~$C([0,T];L^2)$. We consider, for almost
all~$y\in\Omega$, the following linear ordinary differential equation
on~$\pd{f_n}y$
\begin{align} \label{atcouecd} \frac{\partial}{\partial
    t}\brk{\pd{f_n}{y}}&= A \pd{f_n}{y} + B,
\end{align}
where we have introduced the functions
\begin{align}
  A&=\xi|\tau_n| f_{n-1}-f_{n-1} -2\nu f_{n-1} f_n
  &\in C([0,T];L^\infty),\label{rega}\\
  B&=\xi \pd{|\tau_n|}{y} f_{n-1} f_n + \brk{\xi|\tau_n|
    -1}f_n\pd{f_{n-1}}{y} -\nu f_n^2\pd{f_{n-1}}{y} &\in
  C([0,T];L^2).\label{regb}
\end{align}
The Cauchy-Lipschitz Theorem then guarantees the existence
of~$\pd{f_n}{y}(\cdot,y)$ continuous in time, for almost
all~$y\in\Omega$. The Duhamel formula applied to~\eqref{atcouecd}
yields, for all~$t\in[0,T]$ and almost all~$y\in\Omega$,
\begin{align*}
  \pd{f_n}{y}(t,y) = \pd{f_0}{y}(y) e^{\int_0^t A(.,y)} + \int_0^t
  B(s,y) e^{\int_s^t A(.,y)} ds,
\end{align*}
so that, using~\eqref{as:f0nn},~\eqref{rega}
and~\eqref{regb},~$\pd{f_n}{y}$ belong to~$C([0,T];L^2)$.
As~\eqref{atcouecd} is the derivative with respect to~$y$
of~\eqref{eq:atcouec}, this yields~$f_n\in C([0,T];H^1)$.

Now that we have established, for all~$n$, the existence of a
solution~$(u_n,\tau_n,f_n)$ to~\eqref{eq:atcoue} in the appropriate
functional spaces (as in~\eqref{reg}-\eqref{regdt}), we derive,
for~$(u_n,\tau_n,f_n)$, the a priori estimates formally established
on~$(u,\tau,f)$ in the previous steps.  Estimate~\eqref{en1} now reads
\begin{align}\label{aen1} \f12\f{d}{d t}\brk{G\rho
    \Normt{u_n}_\Ld^2+\lambda\Normt{\tau_n}_\Ld^2}+G\eta\Normt{\pd{u_n}{y}}_\Ld^2+
  \Normt{\brk{\sqrt{f_{n-1}}\tau_n}}_\Ld^2
  \nonumber\\
  =G\,a\,\overline{\tau_n}(t).\end{align} Likewise,~\eqref{en2} is now
replaced by
\begin{align}
  \f{d}{d t}\Normt{f_n}_{L^1}&+ \into \brk{f_{n-1} f_n}(t,\cdot) +
  \f\nu2 \into \brk{f_{n-1} f_n^2}(t,\cdot) \leq \frac{\xi^2}{2\nu}
  \Normt{\brk{\sqrt{f_{n-1}}\tau_n}}_\Ld^2,
  \label{aen2}\end{align}
Collecting~\eqref{aen1} and~\eqref{aen2} yields the following estimate,
analogous to~\eqref{ebb},
\begin{align} \f12\f{d}{d t}&\brk{G\rho
    \Normt{u_n}_\Ld^2+\lambda\Normt{\tau_n}_\Ld^2+\frac{2\nu}{\xi^2}
    \Normt{f_n}_{L^1} } \nonumber\\
  &+G\eta\Normt{\pd{u_n}{y}}_\Ld^2+
  \f12\Normt{\brk{\sqrt{f_{n-1}}\tau_n}}_\Ld^2 \le C\,a
  \Normt{\tau_n}_\Ld^2,\label{aebb}
\end{align}
and therefore,~\eqref{eb} holds with~$(u_n,\tau_n,f_n)$ instead of~$(u,\tau,f)$.\\
The arguments given in Step~\ref{s:apf} to derive~\eqref{en3} and in
Step~\ref{s:li} for the~$L^\infty$ estimates can be mimicked for the
approximate system in~$(u_n,\tau_n,f_{n-1})$ instead of~$(u,\tau,f)$,
and the corresponding auxiliary functions~$g_n$ and~$U_n$.

At this point, we have rigorously established on~$(u_n,\tau_n,f_n)$ and
our formal estimates of steps~\ref{s:fee} to~\ref{s:li}:
\begin{align} \sup_n \sup_{t\in[0,T]}
  \brk{\Normt{u_n}_\Ld+\Normt{\tau_n}_\Ld+\Normt{f_n}_{L^1}}\le C_{0,T},
\end{align}
and
\begin{align} \sup_n \sup_{t\in
    [0,T]}\brk{\Normt{U_n}_{H^1}+\Normt{\tau_n}_\Li+\Normt{f_n}_\Li}+\Norm{U_n}_{L^2_T(H^2)}\le
  C_{0,T},
  \label{eq:linfn}\end{align}
where we recall that~$C_{0,T}$ denotes various constants which depend on
the coefficients in system~\eqref{tcoue}, the initial data
$u_0,\tau_0,f_0$ and the time~$T$.

We now turn to the a priori estimates of Step~\ref{s:see}. Using
arguments similar to the formal arguments of Step~\ref{s:see}, we obtain
\begin{align} \label{eq:sapen}
  \f{d}{d t}&\brk{\lambda\Normt{\pd{\tau_n}{y}}_\Ld^2+\Normt{\pd{f_n}{y}}_\Ld^2} \nonumber\\
  &\leq C_{0,T} \brk{
    \Normt{\pd{\tau_n}{y}}_\Ld^2+\Normt{\pd{f_{n-1}}{y}}_\Ld^2+\Normt{\pd{f_n}{y}}_\Ld^2
    +\Normt{\pdd{ U_n}{y}}_\Ld^2}.
\end{align}
We now observe that showing $H^1$ bounds on~$\tau_n$ and~$f_n$ is less
straightforward than in our formal Step~\ref{s:see}. We introduce
\begin{align*}
  Y_n(t)=
  \Norm{\pd{\tau_n}{y}(t,\cdot)}_\Ld^2+\Norm{\pd{f_n}{y}(t,\cdot)}_\Ld^2.
\end{align*}
Integrating~\eqref{eq:sapen} from~$0$ to~$t\le T$, we see that~$Y_n$
satisfies
\begin{align}\label{eq:asape}
  Y_n(t)\le C_{0,T} \int_0^t Y_n + C_{0,T} \int_0^tY_{n-1} + C_{0,T}
  \Norm{U_n}_{L^2_T(H^2)}^2 + Y_0 .
\end{align}
Applying the Gronwall Lemma to~\eqref{eq:asape}, we find
\begin{align*}
  Y_n(t) &\leq \brk{C_{0,T} \Norm{U_n}_{L^2_T(H^2)}^2+ Y_0}~e^{C_{0,T}
    t} + C_{0,T} \int_0^t Y_{n-1}(s) e^{C_{0,T} (t-s)} ds
\end{align*}
which we rewrite
\begin{align*}
  Y_n(t) &\leq C_{0,T} + C_{0,T} \int_0^t Y_{n-1}(s) ds.
\end{align*}
Arguing by induction, one can check that this implies, for
all~$t\in[0,T]$ and~$n$,
\begin{align*}
  Y_n(t)\leq C_{0,T} \sum_{i=0}^{n-1} \frac{(C_{0,T} t)^i}{i!} +
  \frac{(C_{0,T}t)^n}{n!} Y_0.
\end{align*}
It follows that,~$C_{0,T}$ denoting various constants,
\begin{align}\label{eq:sapeyn}
  \sup_n \sup_{t\in[0,T]} Y_n(t)\leq C_{0,T} e^{C_{0,T}T},
\end{align}
Recalling that~$\dfrac{\partial u_n}{\partial y}=\dfrac{\partial
  U_n}{\partial y}-\dfrac1{\eta}(\tau_n-\overline{\tau_n})$, we use
inequalities~\eqref{eq:linfn} and~\eqref{eq:sapeyn} to derive
\begin{align}\label{une}
  \sup_n \sup_{s\in [0,T]}
  \brk{\Norm{u_n(s,\cdot)}_{H^1}+\Norm{\tau_n(s,\cdot)}_{H^1}
    +\Norm{f_n(s,\cdot)}_{H^1}}+\Norm{u_n}_{L^2_{T}(H^2)}\le C_{0,T}.
\end{align}
This implies
\begin{align}\label{uned}
  \sup_n \sup_{s\in [0,T]} \Norm{\pd{u_n}t(s,\cdot)}_{L^2_{T}(L^2)} +
  \Norm{\pd{\tau_n}t(s,\cdot)}_{L^2_{T}(L^2)} +
  \Norm{\pd{f_n}t(s,\cdot)}_{L^2_{T}(L^2)} \le C_{0,T}.
\end{align}

\newstep{Convergence of the sequence of approximate solution}{s:conv}
The bounds obtained in the previous steps, namely~\eqref{une}
and~\eqref{uned} show that, at least up to extraction of a subsequence,
we have the weak convergences
\begin{align}
  (u_{\np},\tau_{\np},f_{\np})\rightharpoonup (u,\tau,f) &\mbox{
    weakly-}\star
  \mbox{  in } L^\infty([0,T];H^1)^3,\label{convf1}\\
  u_{\np}\rightharpoonup u &\mbox{ weakly in } L^2([0,T];H^2),\label{convf2}\\
  \brk{\pd{u_{\np}}t,\pd{\tau_{\np}}t,\pd{f_{\np}}t}\rightharpoonup
  \brk{\pd ut ,\pd \tau t, \pd f t} &\mbox{ weakly in }
  L^2([0,T];L^2)^3.\label{convf3}
\end{align}
But, in order to pass to the limit in~\eqref{eq:atcoue}, we need the
convergence of the whole sequence itself because~\eqref{eq:atcoue}
involves indices~$n-1$ and~$n$ and strong convergence to establish
convergence of the product terms~$f_{n-1} f_n,\ |\tau_n| f_{n-1} f_n,\
f_{n-1} f_n^2$.

We now establish strong convergence of the whole sequence. We prove this
convergence in $\brk{L^\infty([0,T];L^2(\Omega))}^3$. This will a
posteriori imply that all the
convergences~\eqref{convf1},~\eqref{convf2} and~\eqref{convf3} actually
hold true not only for a subsequence, but the whole sequence itself. And
this will provide sufficient information to pass to the limit in our
nonlinear terms.

We introduce the notation:~$\tt{h}_n=h_n-h_{n-1}$ and derive the
evolution equations for $(\tt{u}_n, \tt{\tau}_n, \tt{f}_n)$
\begin{subequations}
  \begin{empheq}[left=\empheqlbrace]{align}
    \rho \frac{\partial \tt{u}_n}{\partial t} &= \eta \frac{\partial^2
      \tt{u}_n}{\partial y^2}
    +   \frac{\partial \tt{\tau}_n}{\partial y}\label{eq:datcoueda},\\
    \lambda \frac{\partial \tt{\tau}_n}{\partial t} &= G \frac{\partial
      \tt{u}_n}{\partial y}
    - f_{n-1}\tt{\tau}_n-\tau_{n-1}\tt{f}_{n-1}  \label{eq:datcouedb},\\
    \frac{\partial \tt{f}_n}{\partial t} & =(-1+\xi
    |\tau_{n-1}|)(f_{n-1}\tt{f}_n+f_{n-1}\tt{f}_{n-1})\nonumber
    \\&\qquad\qquad- \nu f_{n-1}(f_n+f_{n-1})\tt{f}_n-\nu f_{n-1}^2
    \tt{f}_{n-1}+\xi \tt{|\tau_n|}f_{n-1}f_{n} \label{eq:datcouedc}.
  \end{empheq}
  \label{eq:datcoued}
\end{subequations}
Since~$(u_n, \tau_n, f_n)$ belong to the spaces that appear
in~\eqref{reg} and~\eqref{regdt}, the same holds for~$(\tt{u}_n,
\tt{\tau}_n, \tt{f}_n)$.  We multiply
equations~\eqref{eq:datcoueda},~\eqref{eq:datcouedb} and
\eqref{eq:datcouedc}, respectively by~$\tt{u}_n, \tt{\tau}_n$
and~$\tt{f}_n$, integrate over~$\Omega$, sum up and use the
non-negativity of~$f_{n-1}$ and~$f_n$ to find
\begin{align*}
  \f{d}{dt}&\brk{G\rho \Normt{ \tt{u}_n}^2_\Ld + \lambda
    \Normt{\tt{\tau}_n}^2_\Ld +\Normt{\tt{f}_n}^2_\Ld}
  \le -\into \tau_{n-1}\tt{f}_{n-1} \tt{\tau}_n (t,\cdot)\nonumber\\
  &+\into (-1+\xi |\tau_{n-1}|)
  (f_{n-1}\tt{f}_n^2+f_{n-1}\tt{f}_{n-1}\tt{f}_n ) -\nu f_{n-1}^2
  \tt{f}_{n-1} \tt{f}_n +\xi f_{n-1}f_{n} \tt{|\tau_n|} \tt{f}_n
  (t,\cdot).
\end{align*}
The presence of two indices~$n-1$ and~$n$ again makes an additional step
necessary.  We introduce~$X_n(t)=\Normt{\tt{u}_n}^2_\Ld+\Normt{
  \tt{\tau}_n}^2_\Ld+\Normt{\tt{f}_n}^2_\Ld$.  Repeatedly using
the~$L^\infty$-bounds~\eqref{eq:linfn}
on~$\BRK{\tau_n,f_n,\tau_{n-1},f_{n-1}}$ and the Young inequality, we
see that~$X_n$ satisfies
\begin{align}\label{deb}
  \dot{X_n}(t)\le C_{0,T} (X_n(t)+X_{n-1}(t)).
\end{align}
Applying the Gronwall Lemma to~\eqref{deb}, we find \beno X_n(t)\le
C_{0,T}\int_0^t X_{n-1}(s)e^{C_{0,T}(t-s)}ds \le C_{0,T}e^{C_{0,T} T}
\int_0^{t} X_{n-1}(s)ds,\eeno which implies that \beno X_n(t)\le
\f{(C_{0,T} e^{C_{0,T} T}t)^{n-1}}{(n-1)!}  \sup_{s\in[0,T]} X_1(s).
\eeno The sequence~$(u_n, \tau_n, f_n)$ is therefore a Cauchy sequence
in~$\brk{L^\infty([0,T];L^2(\Omega))}^3$. The sequence converges in this
space.

Now that we have strong convergence of the whole sequence, we show how
to pass to the limit in all the terms of~\eqref{eq:atcoue}, including
the nonlinear ones. We only consider~$|\tau_{\np}| f_{\np-1}
f_{\np}$. The other terms can be treated using similar arguments. We use
a classical compactness result~\cite[Theorem~5.1,p. 58]{Lions1969} to
deduce from~\eqref{convf1} and~\eqref{convf3} that~$\tau_{\np}$
and~$f_{\np}$ strongly converge respectively to~$\tau$ and~$f$
in~$L^2([0,T];L^4)^3$. Moreover,~$f_{\np-1}$ strongly converges to~$f$
in~$L^\infty([0,T];L^2)$. We thus have convergence for~$|\tau_{\np}|
f_{\np-1} f_{\np}$ in~$L^1([0,T];L^1)$.

The triple~$(u, \tau, f)$ thus satisfies system~\eqref{eq:atcoue}, at
least in the weak sense. We now derive further regularity. We have
\begin{align*}
  u \in L^2([0,T];H^2) &\mbox{ with } \pd ut \in L^2([0,T];L^2),
\end{align*}
and therefore, by interpolation (see~\cite[Chapter~3,
Lemma~1.2]{Temam1979}),
\begin{align*}
  u \in C([0,T];H^1)\cap L^2([0,T];H^2).
\end{align*}
Moreover, we have
\begin{align*}
  \brk{\pd \tau t, \pd f t} \in L^2([0,T];L^2)^2
\end{align*}
and, using the second a priori estimate~\eqref{eq:sape} ,
\begin{align*}
  \brk{\frac{\partial}{\partial t} \pd \tau y, \frac{\partial}{\partial
      t} \pd fy}\in L^2([0,T];L^2)^2,
\end{align*}
so that,
\begin{align*}
  \brk{\tau,f} \in C([0,T];H^1)^2.
\end{align*}
We have obtained~\eqref{reg} and therefore~\eqref{regdt}, using
system~\eqref{tcoue}.  The non-negativity of the fluidity is preserved,
passing to the limit. This completes the existence proof.

\newstep{Uniqueness}{s:uniq} Consider~$(u_1,\tau_1,f_1)$
and~$(u_2,\tau_2,f_2)$ satisfying~\eqref{reg} and solutions to
system~\eqref{tcoue} supplied with the same initial
condition~$(u_0,\tau_0,f_0)\in H^1(\Omega)$.  We
introduce~$(\tt{u}=u_2-u_1,\tt\tau=\tau_2-\tau_1,\tt{f}=f_2-f_1)$ which
therefore satisfies
\begin{subequations}
  \begin{empheq}[left=\empheqlbrace]{align}
    & \rho \frac{\partial \tt{u}}{\partial t} = \eta \frac{\partial^2
      \tt{u}}{\partial y^2}
    +   \frac{\partial \tt{\tau}}{\partial y}\label{eq:dtcoueda},\\
    & \lambda \frac{\partial \tt{\tau}}{\partial t} = G \frac{\partial
      \tt{u}}{\partial y}
    - f_2\tt{\tau}-\tau_1\tt{f}  \label{eq:dtcouedb},\\
    & \frac{\partial \tt{f}}{\partial t} =-(f_1+f_2)\tt{f} + \xi f_1^2
    \tt{|\tau|} + \xi |\tau_2| (f_1+f_2)\tt{f} -\nu (f_1^2+f_1 f_2 +
    f_2^2) \tt{f}\label{eq:dtcouedc},
  \end{empheq}
  \label{eq:dtcoued}
\end{subequations}
supplied with \emph{homogeneous} boundary conditions and~$(0,0,0)$ as
initial data.  Multiplying
equations~\eqref{eq:dtcoueda},~\eqref{eq:dtcouedb} and
\eqref{eq:dtcouedc}, respectively by~$\tt{u},\ \tt{\tau}$ and~$\tt f$,
integrating over~$\Omega$, summing up, using the~$L^\infty$-bounds
established in Step~\ref{s:li} for terms involving~$\tau_1,\tau_2,f_1,
f_2$ and repeatedly applying the Young inequality, we find
\begin{align*}
  \f12\f{d}{dt}\brk{\rho G\Normt{\tt{u}}_\Ld^2 +
    \lambda\Normt{\tt{\tau}}_\Ld^2 + \Normt{\tt f}_\Ld^2 }\leq C_{0,T}
  \brk{\Normt{\tt{\tau}}_\Ld^2 + \Normt{\tt{f}}_\Ld^2}.
\end{align*}
The Gronwall Lemma then implies uniqueness. This concludes the proof of
Theorem~\ref{existence}.
\begin{flushright}
  \endproof
\end{flushright}
\setcounter{step}{0}
\section{Longtime behaviour for \emph{homogeneous} boundary conditions}
\label{sec:hom}
In this section, we study the longtime behaviour of
system~\eqref{eq:couette} supplied with \emph{homogeneous} boundary
conditions. We will show convergence to a steady state and establish a
rate for this convergence. For \emph{homogeneous} boundary conditions,
the $H^1$-steady states of~\eqref{eq:couette} such that~$f\ge0$ are
exactly the states~$(u\equiv 0,\tau\equiv c, f\equiv 0),$ where~$c$ is a
constant throughout the domain.

Indeed, such a steady state~$(\uinf,\tauinf,\finf)$ satisfies, combining
equation~\eqref{eq:couettea} integrated over the domain
and~\eqref{eq:couetteb},
\begin{align} \label{eq:c} \tauinf \brk{\frac\eta G \finf+ 1}= c,
\end{align}
where $c$ is a constant over the domain. We now distinguish between two
cases. Either~$c=0$, in which case~$\tauinf\equiv0$. The
\emph{homogeneous} boundary conditions on~$u$ and~\eqref{eq:couetteb}
imply that~$\uinf\equiv0$ and~\eqref{eq:couettec}
that~$\finf\equiv0$. Or~$c\neq 0$ and it follows from~\eqref{eq:c}
that~$\tauinf$ is non-zero and has a constant sign and
from~\eqref{eq:couetteb} that~$\pd \uinf y$ has a constant sign.
Because of the \emph{homogeneous} boundary conditions on the velocity,
we obtain that~$\uinf \equiv 0$. Therefore,~\eqref{eq:couetteb}
yields~$\finf\tauinf=0$ and~$\finf\equiv 0$, because~$\tauinf$ is
non-zero in this case.

We will show that the longtime behaviour differs both in terms of steady
state and rate of convergence, depending whether~$f_0 \nequiv 0$ or~$f_0
\equiv 0$. When~$f_0 \nequiv 0$, a case studied in
subsection~\ref{sec:longhom1}, the solution~$(u,\tau,f)$ converges to
the steady state~$(0,0,0)$ in the longtime and the rates of convergence
are power-laws of the time. In the case~$f_0 \equiv 0$, the fluidity~$f$
vanishes for all time, as easily seen on~\eqref{eq:couettec}.  In
subsection~\ref{sec:longhom0}, we show that~$(u,\tau,f)$ then converges
to~$(0,\overline{\tau_0},0)$ in the longtime at an exponential
rate,~$\overline{\tau_0}$ being the average of~$\tau_0$ over~$\Omega$.
Evidently, the former case~$f_0 \nequiv 0$ require more efforts than the
latter case~$f_0 \equiv 0$ where~$f\equiv 0$ for all times.

\subsection{Case~$f_0 \nequiv 0$}
\label{sec:longhom1}
In this subsection, we consider the case~$f_0\ge 0,\ f_0 \nequiv 0$. We
first establish the convergence in the longtime.

\begin{thm} \label{th:longhom1} Supply system~\eqref{eq:couette} with
  homogeneous boundary conditions and initial conditions that
  satisfy~\eqref{as:f0nn} and~$f_0 \nequiv 0$.
  The solution $(u,\tau,f)$, the existence and uniqueness of which have
  been established in Theorem~\ref{existence}, converges to the steady
  state~$(0,0,0)$ in~$H^1(\Omega)\times L^\infty(\Omega) \times
  L^\infty(\Omega)$ in the longtime:
  \begin{align*}
    \Normt{u}_{H^1} +\Normt{\tau}_\Li +\Normt{f}_\Li\rightarrow 0 .
  \end{align*}
\end{thm}

\proof The proof falls in three steps. In the first step, we establish a
lower bound for the average of the fluidity~$f$, which, in
Step~\ref{sec:l2stab2}, is useful to prove convergence in the longtime
in~$L^2(\Omega)$.  In the third step, we show convergence
of~$(u(t,\cdot),\tau(t,\cdot),f(t,\cdot))$ in~$H^1(\Omega)\times
L^\infty(\Omega)\times L^\infty(\Omega)$.

In this section,~$C_0$ denotes various constants that are independent
from time, while~$C_i$, $i$=1,...,4 denote some fixed constants
independent from time. These constants~$C_0$ and~$C_i$ used to be
denoted~$C_{0,T}$ in the previous section. The subscript~$T$ is now
omitted because, as explained in Remarks~\ref{rmk:chom1}
and~\ref{rmk:chom2}, the constants are independent from~$T$
for~\emph{homogeneous} boundary conditions.

\newstep{A lower bound for the average of~$f$}{s:low} We first derive a
lower bound on~$\bar{f}$, defined as in~\eqref{eq:avg}, and not directly
on~$f$ because the latter may vanish (since~$f_0$ may vanish) on some
part of the domain.  Since~$f_0 \nequiv 0$, there exists, by continuity
of~$f_0$ (assumed in~$H^1$), a non-empty closed interval~$\Omega_0$
in~$\Omega$ where~$f_0$ does not vanish.  Arguing as in Step~\ref{s:pos}
of the proof of Theorem~\ref{existence}, the fluidity~$f$ does not
vanish for all~$t>0$ and~$y\in\Omega_0$.  The evolution
equation~\eqref{eq:couettec} on~$f$ rewrites, for all~$t>0$
and~$y\in\Omega_0$,
\begin{align}
  \label{equ-f} \frac{\partial }{\partial t} \f1{f}=1-\xi|\tau|+\nu f.
\end{align}
As explained in Remark~\ref{rmk:chom2}, the~$L^\infty$-bounds on~$\tau$
and~$f$ are uniform in time for \emph{homogeneous} boundary conditions.
The equation~\eqref{equ-f} thus implies, for all~$y\in\Omega_0$
and~$t>0$, \beno \frac{\partial }{\partial t} \f1{f}\le C_0,\eeno and
therefore,
\begin{align}\label{eq:lf}
  f(t,y)&\ge \f{1}{\f{1}{f_0(y)}+C_0 t},  \nonumber\\
  & \ge\f{1}{\Norm{\f{1}{f_0}}_{\Li(\Omega_0)}+C_0 t}.
\end{align}
Since~$ \bar{f}\ge \dps \int_{\Omega_0}f$, this yields the lower bound
\begin{align}
  \bar{f}(t)\ge \frac{C_1}{1+C_0t}.
\end{align}

\newstep{Longtime convergence in~$L^2(\Omega) \times L^2(\Omega) \times
  L^2(\Omega)$}{sec:l2stab2} We now show the longtime convergence
in~$L^2$. Estimates~\eqref{en1} and~\eqref{en3} respectively rewrite,
for \emph{homogeneous} boundary conditions,
\begin{align}\label{eq:lot1} \f12\f{d}{d t}\brk{G\rho
    \Normt{u}_{L^2}^2+\lambda\Normt{\tau}_{L^2}^2}+G\eta\Normt{\pd{u}{y}}_{L^2}^2+
  \Normtp{\sqrt{f}\tau}_{L^2}^2=0\end{align} and
\begin{align}\label{eq:lot2}
  \f12 \f{d}{d t}\Normt{\pd{U}{y}}_\Ld^2+
  \f{\eta}{2\rho}\Normt{\pdd{U}{y}}_\Ld^2\le C_2\brk{\Normt{f}_{L^1}
    \Normtp{\sqrt{f}\tau}_\Ld^2+\Normt{u}_\Ld^2},
\end{align}
where~$U$ is defined by~\eqref{def:U}.  The evolution equation
on~$\bar{\tau}$ writes
\begin{align} \label{eq:bartau} \lambda
  \deriv{\bar{\tau}}{t}+\bar{f}\bar{\tau}=-\overline{f(\tau-\bar{\tau})}.
\end{align}
We introduce the positive scalar~$\epsilon$, to be fixed later on.  We
use the Cauchy-Schwartz and Young inequalities
\begin{align*} \Abs{\bar\tau} \Abs{\overline{f(\tau-\bar{\tau})}} \le
  \sqrt{\bar{f}}\Abs{\bar\tau}\Normtp{\sqrt{f}(\tau-\bar{\tau})}_\Ld \le
  \epsilon \bar{f} \Abs{\bar\tau}^2 + \f1{4\epsilon}
  \Normtp{\sqrt{f}(\tau-\bar{\tau})}_\Ld^2 ,
\end{align*}
so that, multiplying evolution equation~\eqref{eq:bartau} by~$\bar\tau$,
we obtain
\begin{align}\label{eq:lot4}\f\lambda2\deriv{|\bar{\tau}|^2}{t}(t)+(1-\epsilon)
  \bar{f}|\bar{\tau}|^2(t)\le
  \f1{4\epsilon}\Normt{f}_\Li\Normtp{\tau-\bar{\tau}}_\Ld^2.
\end{align}
The evolution equation on~$\tau-\bar{\tau}$ reads
\begin{align} \label{eq:tmt} \lambda \frac{\partial }{\partial
    t}(\tau-\bar{\tau})+\f{G}{\eta}(\tau-\bar{\tau})=-(f\tau-\overline{f\tau})+G\pd{U}{y}.
\end{align}
Multiplying evolution equation~\eqref{eq:tmt} by~$\tau-\bar{\tau}$,
integrating over~$\Omega$ and repeatedly using the Young inequality, we
find
\begin{align*}
  \f12\f{d}{d
    t}(\lambda\Normtp{\tau-\bar{\tau}}_\Ld^2)+\f{G}{2\eta}\Normtp{\tau-\bar{\tau}}_\Ld^2\le
  C\brk{\Normtp{f\tau}_\Ld^2 + \Normt{\pd{U}{y}}_\Ld^2},
\end{align*}
so that, using the uniform in time~$\Li$-bound on~$f$,
\begin{align}\label{eq:lot3} \f12\f{d}{d
    t}(\lambda\Normtp{\tau-\bar{\tau}}_\Ld^2)&+\f{G}{2\eta}\Normtp{\tau-\bar{\tau}}_\Ld^2  \nonumber\\
  &\le
  C_3\brk{\Normt{f}_\Li\Normtp{\sqrt{f}\tau}_\Ld^2+\Normt{\pd{U}{y}}_\Ld^2
  }. \end{align} We introduce some positive scalars~$\mm_1, \mm_2,
\mm_3$ and the energy function
\begin{align} \label{def:E} E(t)=& \mm_1 (G\rho\Normt{u}_\Ld^2+\lambda
  \Normt{\tau}_\Ld^2)+\mm_2 \Normt{\pd{U}{y}}_\Ld^2 \nonumber\\& + \mm_3
  \lambda \Normtp{\tau-\bar{\tau}}^2_\Ld+\lambda |\bar{\tau}(t) |^2,
\end{align}
which therefore satisfy, combining \eqref{eq:lot1}, \eqref{eq:lot2},
\eqref{eq:lot4} and~\eqref{eq:lot3},
\begin{align}\label{opene} \f12\f{dE}{d t}(t)
  &+\brk{\mm_1-C_2\mm_2\Normt{f}_{L^1}-C_3\mm_3\Normt{f}_\Li}\Normtp{\sqrt{f}\tau}_\Ld^2\nonumber
  \\&+ (\mm_1G\eta-C_2\mm_2C_p)\Normt{\pd{u}{y}}_\Ld^2+
  \brk{\f{\eta}{2\rho}\mm_2-C_3\mm_3C_p}\Normt{\pdd{U}{y}}_\Ld^2\nonumber\\
  &+\brk{\f{G}{2\eta}\mm_3-\f1{4\epsilon}\Normt{f}_\Li}\Normtp{\tau-\bar{\tau}}^2_\Ld+(1-\epsilon)
  \bar{f}|\bar{\tau}|^2(t)\le 0, \end{align} where~$C_p$ is the
Poincar\'e constant.

The coefficients~$\mm_1, \mm_2, \mm_3$ are chosen sufficiently large so
that, for all time~$t>0,$ every term in the left-hand side
of~\eqref{opene} is positive. The conditions
\begin{align*}
  \f{G}{2\eta}\mm_3&>\f1{4\epsilon}\sup_{t>0}\Normt{f}_\Li,\\
  \f{\eta}{2\rho}\mm_2&>C_3\mm_3C_p,\\
  \mm_1&>\max\brk{\frac{C_2\mm_2C_p}{G\eta},
    C_2\mm_2\sup_{t>0}\Normt{f}_{L^1}+C_3\mm_3\sup_{t>0}\Normt{f}_\Li}
\end{align*}
are sufficient.  Using in addition the lower bound~\eqref{eq:lf} and the
Poincar\'e inequality,~\eqref{opene} becomes
\begin{align*}
  \f12
  \f{dE}{dt}(t)+C_0\brk{\Normt{u}_\Ld^2+\Normtp{\tau-\bar{\tau}}^2_\Ld+\Normt{\pd{U}{y}}_\Ld^2}+
  (1-\epsilon) \f{C_1}{1+C_0 t} |\bar{\tau}|^2(t)\le 0.
\end{align*}
Using the triangle inequality
\begin{align}\label{eq:tri}
  \frac{\epsilon}{2\mm_1} \mm_1 \Normt{\tau}_\Ld^2 \le \epsilon
  \Abs{\bar{\tau}}^2(t)+ \epsilon \Normtp{\tau-\bar{\tau}}_\Ld^2,
\end{align}
we find, for~$t$ sufficiently large,
\begin{align}\label{lot6}
  \f12 \f{dE}{d t} + \f1\lambda
  \min\brk{1-2\epsilon,\f{\epsilon}{2\mm_1}} \f{C_1}{1+C_0 t}E\le 0.
\end{align}
We take~$\epsilon<\dfrac12$ and apply the Gronwall Lemma to~\eqref{lot6}
to obtain that~$E$ goes to zero in the longtime limit. In particular, we
have
\begin{align} \label{de1} \lim_{t\rightarrow\infty}
  \brk{\Normt{u}_{L^2}^2+\Normt{\tau}_{L^2}^2+\Normt{\pd{U}{y}}_{L^2}^2}=0.
\end{align}

\newstep{Longtime convergence in~$H^1(\Omega)\times
  L^\infty(\Omega)\times L^\infty(\Omega)$}{sec:linfstab}
Combining~\eqref{en3} and~\eqref{en4}, using the uniform in
time~$\Li$-bound on~$f$ and the Poincar\'e inequality on~$\pd Uy$, the
spatial average of which is zero, we obtain
\begin{align*}
  \f{d}{d t}&\brk{\Normt{\pd{U}{y}}_{L^2}^2+\lambda
    |\tau(t,y)|^2}+ \brk{\Normt{\pdd{U}{y}}_{L^2}^2+|\tau(t,y)|^2}\\
  &\leq C_0
  \brk{\Normt{u}_{L^2}^2+\Normt{\tau}_{L^2}^2+\Normt{\pd{U}{y}}_{L^2}^2}.
\end{align*}
We apply the Gronwall Lemma and use the convergence~\eqref{de1} that is
uniform in space to derive
\begin{align*}
  \lim_{t\rightarrow\infty} \Normt{\tau}_\Li=0.
\end{align*}
Using the convergence of~$\Normt{\tau}_\Li$, the evolution
equation~\eqref{eq:couettec} on~$f$ implies, for~$t$ sufficiently large,
\begin{align}
  \pd f t \le -\f12 f^2.
\end{align}
This yields the convergence of~$\Normt{f}_\Li$ to zero in the longtime.

Additionally, using the definition~\eqref{def:U} and~\eqref{de1}, $\pd u
y$ converges to zero in~$L^2(\Omega)$. This ends the proof.
\setcounter{step}{0}
\endproof

We now turn to making precise the rates of convergence to the
steady-state.  We introduce the non-negative scalar
\begin{align} \label{def:beta} \beta=\mbox{meas}\BRK{y\in \Omega
    |f_0(y)>0}.
\end{align}
By assumption in this section, we have~$\beta>0$.  The following result
establishes the convergence rates in function of~$\beta$. In
Section~\ref{ssec:numh}, we will check using numerical simulations that
these rates are indeed sharp.

\begin{thm} \label{th:longhomt} Supply system~\eqref{eq:couette} with
  homogeneous boundary conditions and initial conditions that
  satisfy~\eqref{as:f0nn} and~$f_0\nequiv 0$.  The solution~$(u, \tau,
  f)$, the existence and uniqueness of which have been established in
  Theorem~\ref{existence}, satisfies the following convergence
  estimates: for any arbitrarily small~$\vveps>0$, there exists a
  constant~$\kgp_\vveps$ independent from time and there exists a
  time~$t_0$, both depending on the domain, the initial data, the
  coefficients in the system and~$\vveps$, such that, for all~$t>t_0$,
  \begin{align}\label{dr2}
    \Normt{u}_{H^1}+\Normt{\tau}_\Ld &\le \kgp_\vveps
    (1+t)^{-\f{\beta}{\lambda}(1-\vveps)},
  \end{align}
  where~$\beta$ is defined by~\eqref{def:beta} and for all~$t>t_0$ and
  $y\in\Omega~$, we have
  \begin{align} \label{vvdr3} \f1{\f1{f(t_0,y)}+(1+\vveps)(t-t_0)}\le
    f(t,y)\le \f{1}{\f1{f(t_0,y)}+ (1-\vveps)(t-t_0) }.
  \end{align}
  In addition, there exists another constant~$\kgp_\vveps$, such that,
  for all~$t>t_0$,
  \begin{align}
    \Normt{u}_{H^1}+\Normtp{\tau-\overline{\tau}}_\Ld&\le
    \kgp_\vveps (1+t)^{-1-\f{\beta}{\lambda}(1-\vveps)},\label{dr1}\\
    \Normtp{\eta \pd{u}{y}+\tau-\overline{\tau}}_\Ld&\le \kgp_\vveps
    (1+t)^{-2-\f{\beta}{\lambda}(1-\vveps)},\label{dr4}
  \end{align}
  where the function~$\bar{\tau}$
  is the spatial average of~$\tau$, defined as in~\eqref{eq:avg}.\\
\end{thm}

\proof The proof falls in four steps. We first consider the fluidity,
then derive first convergence rates for the velocity and the stress. A
study of the auxiliary function defined by~\eqref{def:U} next allows to
conclude on the convergence estimates~\eqref{dr1} and~\eqref{dr4}.

We fix~$\veps$ an arbitrarily small positive scalar, actually equal
to~$\dfrac\vveps4$, where~$\vveps$ is the constant that appears in the
statement of the Theorem. The constants~$\kgp_\veps$ depend on~$\veps$
and have value that may vary from one instance to another, the actual
value being irrelevant.

\newstep{Convergence rate for the fluidity}{s:cf} In view of
Theorem~\ref{th:longhom1},~$\Normt{\tau}_\Li$ and $\Normt{f}_\Li$ vanish
in the longtime. This implies that there exists a time~$t_0$, such that,
for all~$y\in\Omega$ and~$t>t_0$, the evolution
equation~\eqref{eq:couettec} on~$f$ leads to
\begin{align} \label{ineq:f} -(1+\epsilon) f^2(t,y) \le \pd f t(t,y) \le
  - (1-\epsilon) f^2(t,y).
\end{align}
For all~$y\in\Omega$ such that~$f_0(y)>0$, we have~$f(t_0,y)>0$, as
shown in Step~\ref{s:pos} of the proof of Theorem~\ref{existence}.  The
equation~\eqref{ineq:f} becomes, for such~$y$ and~$t>t_0$,
\begin{align} \label{dr3} \f1{\f1{f(t_0,y)}+(1+\veps)(t-t_0)}\le
  f(t,y)\le \f{1}{\f1{f(t_0,y)}+ (1-\veps)(t-t_0) }.
\end{align}
and hence~\eqref{vvdr3}.  This inequality is also valid for all~$y$ such
that~$f_0(y)=0$, that is~$f(t_0,y)=0$, and therefore for
all~$y\in\Omega$.

\newstep{First convergence rates for the velocity and the stress}{s:cs}
We first make more precise the lower bound on~$\bar{f}$. As~$f_0$ is
continuous, there exists a closed set~$\Omega_\epsilon$ such
that~$f_0(\Omega_\epsilon)> 0$ and~$\mbox{meas}\{\Omega_\epsilon\}
=\beta(1-\epsilon)$. As shown in Step~\ref{s:pos} of the proof of
Theorem~\ref{existence}, we also have~$f(t_0,\Omega_\epsilon)> 0$.
Furthermore, as~$f(t_0,\cdot)$ is continuous, we obtain
$f(t_0,\Omega_\epsilon)>\kgp_\veps$. The inequality~\eqref{dr3} thus
becomes, for all~$y\in\Omega_\veps$ and~$t>t_0$,
\begin{align*}
  f(t,y)>\frac1{1+\veps} \frac1{\kgp_\veps+t}.
\end{align*}
It follows from~$\bar f \ge \dps\int_{\Omega_\epsilon} f$ that, for
all~$t>t_0$,
\begin{align}\label{eq:lowfb}
  \bar{f}(t)\ge \beta \frac{1-\epsilon}{1+\veps}\f1{\kgp_\veps+t} .
\end{align}
We now use the energy~$E$ introduced in~\eqref{def:E}.
As~$\Normt{f}_\Li$ vanishes in the longtime, the coefficients~$\mm_1,
\mm_2, \mm_3$ can be chosen arbitrarily small in~\eqref{opene},
independently from~$\epsilon$, for~$t>t_0$ sufficiently large.  We
insert~\eqref{eq:lowfb} in~\eqref{opene} so that, for sufficiently
large~$t$,
\begin{align*}
  \f12 \f{dE}{dt}(t)&+C_0\brk{\Normt{u}_\Ld^2+
    \Normtp{\tau-\bar{\tau}}^2_\Ld+\Normt{\pd{U}{y}}_\Ld^2} \nonumber\\&
  +(1-\epsilon) \beta\frac{1-\epsilon}{1+\veps} \f1{\kgp_\veps+t}
  |\bar{\tau}|^2(t)\le 0.
\end{align*}
Using the triangle inequality~\eqref{eq:tri}, we obtain, for
sufficiently large~$t$,
\begin{align*}
  \f12 \f{dE}{d t} +\frac{\beta}{\lambda}
  \min\brk{1-2\epsilon,\f{\epsilon}{2\mm_1}}
  \frac{1-\epsilon}{1+\epsilon} \f1{\kgp_\veps+t} E\le 0.
\end{align*}
and therefore, using that~$1-2\epsilon<\dfrac\epsilon{2\mm_1}$
as~$\mm_1$ is arbitrarily small,
\begin{align} \label{eq:rateE} \f12 \f{dE}{d t} +\frac{\beta}{\lambda}
  \brk{1-4\epsilon} \f1{\kgp_\veps+t} E\le 0.
\end{align}
Applying the Gronwall Lemma to~\eqref{eq:rateE}, we find,
\begin{align}\label{de4}
  \Normt{u}_\Ld^2+\Normt{\tau}_\Ld^2+\Normt{\pd{U}{y}}_\Ld^2 \le
  \kgp_\veps (1+t)^{-2\f{\beta}{\lambda}(1-4\veps)} .
\end{align}
where we recall that~$\kgp_\veps$ denotes various constants.  We have
obtained~\eqref{dr2}.

\newstep{Convergence rate for the auxiliary function~$U$}{s:dwe} We
recall that the function~$U$ is defined by~\eqref{def:U}.  We first
prove that $U$ is more regular than claimed in~\eqref{reg:U}. We
rewrite~\eqref{eq:U}
\begin{align} \label{req:U} \pd{U}{t}- \f{\eta}{\rho}
  \pdd{U}{y}=-\f{1}{\lambda\eta}\int_0^y(f\tau-\overline{f\tau})d
  x+\frac{G}{\lambda\eta}u.
\end{align}
We deduce that $\pd Uy$ satisfies the heat equation with a right-hand
side in $\Ldloc((t_0,+\infty),L^2)$ and initial condition $\pd Uy
(t_0,\cdot)\in H^1(\Omega)$ at time $t_0$ (up to a possible modification
on a set of times of measure zero).  Therefore, we have $\pd Uy\in
H^1_{loc}((t_0,+\infty),L^2),$ so that
\begin{align} \label{rreg:U} U\in H^1_{loc}((t_0,+\infty),H^1_0).
\end{align}

We next differentiate~\eqref{eq:U} with respect to~$t$,
insert~\eqref{eq:couettea} and find,
\begin{align} \label{eq:dwe} \pdd{U}{t}- \f{\eta}{\rho} \frac{\partial^2
  }{\partial y^2}\brk{\pd U t}-\frac{G}{\lambda\rho}\pdd{U}{y}= I,
\end{align}
where~$I$ is the function defined by
\begin{align} \label{eq:I} I(t,y)=
  -\f{1}{\lambda\eta}\int_0^y\brk{\pd{f\tau}{t}-
    \pd{\overline{f\tau}}{t}}d x.
\end{align}
We now regularize $I$ as follows.  We consider a sequence of functions
$I_m$ such that for all $m$, $I_m$ is infinitely differentiable from
$(t_0,+\infty)$ to $L^2(\Omega)$ and as $m\rightarrow\infty$,
\begin{align}
  I_m\rightarrow I \mbox{ in } \Ldloc((t_0,+\infty),L^2)
  \label{conv:I_m}\end{align}
Consider a solution $U_m\in C^\infty((t_0,+\infty),H^2\cap H^1_0)$ to
\begin{align} \label{eq:dwem} \pdd{U_m}{t}- \f{\eta}{\rho}
  \frac{\partial^2 }{\partial y^2}\brk{\pd{U_m}
    t}-\frac{G}{\lambda\rho}\pdd{U_m}{y}= I_m.
\end{align}
Equation~\eqref{eq:dwem} has been studied
in~\cite{hx,kawa-shi}. Inspired by arguments from these references, we
introduce the energy functions~$H_m$ and~$F_m$ depending on a
constant~$\delta\in(0,1)$ to be determined later
\begin{align*} H_m(t)
  =\Normt{\pd{U_m}{t}}_\Ld^2+&\frac{G}{\lambda\rho}\Normt{\pd{U_m}{y}}_\Ld^2+
  \f{\eta}{\rho}\delta \Normt{\pd{U_m}{y}}_\Ld^2+ \delta \into
  \brk{\pd{U_m}{t} U_m}(t,\cdot)
  \\&+2\f{\eta}{\rho}\delta^2\Normt{\frac{\partial}{\partial
      t}\pd{U_m}{y}}_\Ld^2+2\frac{G}{\lambda\rho}\delta^2 \into
  \brk{\frac{\partial}{\partial t}\brk{\pd{U_m}{y}}
    \pd{U_m}{y}}(t,\cdot)
\end{align*}
and
\begin{align*}
  F_m(t)=\brk{\f{\eta}{\rho}-\delta^2\frac{G}{\lambda\rho}}
  \Normt{\frac{\partial}{\partial t}\pd{U_m}{y}}_\Ld^2
  &-\delta\Normt{\pd{U_m}{t}}_\Ld^2 \\&
  +\delta\frac{G}{\lambda\rho}\Normt{\pd{U_m}{y}}_\Ld^2
  +\delta^2\Normt{\pdd{U_m}{t}}_\Ld^2.
\end{align*}
We multiply~\eqref{eq:dwem} by $\pd{U_m}{t}+\delta
U_m+\delta^2\pdd{U_m}{t}$, integrate over $\Omega$ and find
\begin{align*}
  \f12\f{dH_m}{d t}(t)+F_m(t)=\into I_m(t,y) \brk{\pd{U_m}{t}+\delta
    U_m+\delta^2\pdd{U_m}{t}}(t,y)dy.
\end{align*}
We use the Poincar\'e inequality and choose~$\delta$ sufficiently small,
depending on the domain and the coefficients in~\eqref{eq:couette} such
that, for suitable constants~$c_1$,~$c_2$ and~$c_3$,
\begin{align} \label{enc:H} c_1\brk{\Normt{\frac{\partial}{\partial t}
      \pd{U_m}{y}}_\Ld^2+\Normt{\pd{U_m}{y}}_\Ld^2} \le H_m(t)\le
  c_2\brk{\Normt{\frac{\partial}{\partial t}
      \pd{U_m}{y}}_\Ld^2+\Normt{\pd{U_m}{y}}_\Ld^2}
\end{align}
and
\begin{align} \label{lb:F} F_m(t)\ge c_3\brk{
    \Normt{\frac{\partial}{\partial t}
      \pd{U_m}{y}}_\Ld^2+\Normt{\pd{U_m}{y}}_\Ld^2 +
    \Normt{\pdd{U_m}{t}}_\Ld^2}.
\end{align}
Using the upper bound in~\eqref{enc:H} and~\eqref{lb:F} and the Young
and the Poincar\'e inequalities, we obtain
\begin{align}
  \f{dH_m}{d t}(t)+ C H_m(t)\le \Normt{I_m}_\Ld^2.
\end{align}
We multiply the above equation by~$e^{Ct}$, integrate from $t_0$ to $t$
and find
\begin{align} \label{int:Hm} H_m(t) e^{Ct} \le H_m(t_0) e^{Ct_0} +
  \int_{t_0}^t \Norm{I_m(s,\cdot)}_\Ld^2 e^{Cs} ds.
\end{align}
Equation~\eqref{eq:dwe} is linear so that by~\eqref{rreg:U}
and~\eqref{conv:I_m}, we can pass to the limit $m\rightarrow\infty$
in~\eqref{int:Hm} and find, for all $t>t_0$,
\begin{align} \label{int:H} H(t) e^{Ct} \le H(t_0) e^{Ct_0} +
  \int_{t_0}^t \Norm{I(s,\cdot)}_\Ld^2 e^{Cs} ds.
\end{align}
where $H$ is defined as $H_m$ with $U$ instead of $U_m$.

The study of~\eqref{eq:dwe} reduces to the understanding
of~\eqref{int:H}.  We now make precise the behaviour of~$I$ or more
precisely at the one of~$\pd{f\tau}t$.  We combine
equations~\eqref{eq:couetteb} and~\eqref{eq:couettec} to find
\begin{align}\label{eq:ftau}
  \pd{f\tau}t=f\brk{-\f1\lambda f\tau+\frac{G}\lambda \pd u y}+\tau
  \brk{-f^2-\nu f^3+\xi |\tau| f^2} .
\end{align}
Multiplying the evolution equation~\eqref{eq:ftau} by~$f\tau$ and
integrating over~$\Omega$ yields
\begin{align}\label{no}
  \f12 \Normt{\pd{f\tau}t}_\Ld^2 &= \into \brk{-\f1\lambda-1-\nu
    f+\xi|\tau|} f^3\tau^2 + \frac{G}\lambda \into f^2\tau \pd u
  y\nonumber\\
  &\le C_0 \Normt{f}_\Li^2 \brk{\Normt{\pd u
      y}_\Ld^2+\Normt{\tau}_\Ld^2},
\end{align}
where we have used the~$L^\infty$-bounds on both~$\tau$ and~$f$ and the
Cauchy-Schwarz inequality to derive the second line.
Inserting~\eqref{dr3} which gives the convergence in~$\dfrac1t$
of~$\Normt{f}_\Li$ and~\eqref{de4}, equation~\eqref{no} implies
\begin{align} \label{taux:ftau} \Normt{\pd{f\tau}t}_\Ld^2\le \kgp_\veps
  (1+t)^{-2-2\f{\beta}\lambda(1-4\veps)}.
\end{align}
Since the~$L^2$-norm of~$\pd{f\tau}t$ controls the~$L^2$-norm of~$I$, we
insert~\eqref{taux:ftau} in~\eqref{int:H} so that, for all~$t>t_0$,
\begin{align} \label{eq:gronH} H(t)e^{Ct} \le H(t_0)e^{C t_0} +
  \kgp_\veps \int_{t_0}^t
  \frac{e^{Cs}}{(1+s)^{2+2\f{\beta}\lambda(1-4\veps)}} ds.
\end{align}
Moreover, for~$q>0$, for all~$t>t_0$, we integrate by parts to obtain
\begin{align} \label{ipp:et/t} \int_{t_0}^t \frac{e^{Cs}}{(1+s)^q}ds \le
  \frac{q}{C(1+t_0)} \int_{t_0}^t \frac{e^{Cs}}{(1+s)^q}ds
  +\frac{e^{Ct}}{C(1+t)^q}.
\end{align}
We insert~\eqref{ipp:et/t} with~$q=2+2\f{\beta}\lambda(1-4\veps)$
in~\eqref{eq:gronH}, so that for $t$ sufficiently large
\begin{align}\label{taux:H}
  H(t) \le \kgp_\veps \frac{1}{(1+t)^{2+2\f{\beta}\lambda(1-4\veps)}}.
\end{align}
Using the lower bound in~\eqref{enc:H}, we have therefore obtained
\begin{align}\label{de5}
  \Normt{\pd U y}_\Ld^2\le \kgp_\veps
  (1+t)^{-2-2\f{\beta}\lambda(1-4\veps)}.
\end{align}

\newstep{Convergence rates~\eqref{dr1} and~\eqref{dr4}}{s:cv}
Using~\eqref{eq:lot3} rewritten as
\begin{align*}
  \f12\f{d}{d
    t}&(\lambda\Normtp{\tau-\bar{\tau}}_\Ld^2)+\f{G}{2\eta}\Normtp{\tau-\bar{\tau}}_\Ld^2
  \\&\le
  C\brk{\Normt{\pd{U}{y}}_\Ld^2+\Normt{f}_\Li^2\Normt{\tau}_\Ld^2},
\end{align*}
and convergence estimates~\eqref{dr3},~\eqref{de4},~\eqref{de5}, we
obtain \beno \Normtp{\tau-\bar{\tau}}_\Ld^2\le \kgp_\veps
(1+t)^{-2-2\f{\beta}\lambda(1-4\veps)}, \eeno and eventually \beno
\Normt{\pd{u}{y}}_\Ld^2\le \kgp_\veps
(1+t)^{-2-2\f{\beta}\lambda(1-4\veps)}.\eeno We thus obtain~\eqref{dr1}
with $\vveps=4\veps$ and conclude establishing~\eqref{dr4} as follows:
we return to~\eqref{no} and improve the convergence estimate for
$\pd{f\tau}t$ , namely \beno \Normt{\pd{f\tau}t}_\Ld^2\le \kgp_\veps
(1+t)^{-4-2\f{\beta}\lambda(1-4\veps)}.  \eeno This implies,
mimicking~\eqref{eq:gronH} and using~\eqref{ipp:et/t} with
$q=4+2\f{\beta}\lambda(1-4\veps)$, that for~$t$ sufficiently large,
\begin{align}
  \Normt{\pd U y}_\Ld^2\le \kgp_\veps
  (1+t)^{-4-2\f{\beta}\lambda(1-4\veps)},
\end{align}
that is~\eqref{dr4} with $\vveps=4\veps$.  \setcounter{step}{0}
\endproof

\subsection{Case~$f_0 \equiv 0$}
\label{sec:longhom0}
In the case~$f_0\equiv 0$,~$f$ vanishes for all
time. System~\eqref{eq:couette} then reads
\begin{subequations}
  \begin{empheq}[left=\empheqlbrace]{align}
    \rho \pd{u}{t} &= \eta \pdd{u}{y}+ \pd{\tau}{y}, \label{tcoue0a}\\
    \lambda \pd{\tau}{t} &= G \pd{u}{y}.\label{tcoue0b}
  \end{empheq}
  \label{od}
\end{subequations}
The existence and uniqueness of a regular solution to~\eqref{od} is easy
to establish.  The longtime behaviour of system~\eqref{od} is now made
precise.

\begin{thm}
  Supply system~\eqref{od} with homogeneous boundary conditions.
  Consider a solution~$(u, \tau)$ in the space
  \begin{align*}
    \brk{C([0,+\infty);H^1)\cap \Ldloc([0,+\infty);H^2)}\times
    C([0,+\infty);H^1)
  \end{align*}
  Then, the solution converges exponentially fast to the steady
  state~$(0,\overline{\tau_0})$ in~$H^1(\Omega)\times L^2(\Omega)$ in
  the longtime: there exist two constants~$C$, independent from time and
  initial data, and~$C_0$, independent from time, such that, for~$t$
  sufficiently large ,
  \begin{align}\label{dr5}
    \Normt{\pd uy}_\Ld+\Norm{\tau(t,\cdot)-\overline{\tau_0}}_\Ld\le C_0
    e^{-C t}.
  \end{align}
\end{thm}

\proof We perform the same manipulations as those used to obtain
equation~\eqref{eq:dwe} in Step~\ref{s:dwe} of the proof of
Theorem~\ref{th:longhomt}. Since we deal here with the case~$f\equiv 0$,
we have~$I=0$ in~\eqref{eq:dwe}.  We have proven that studying the
longtime behaviour to~\eqref{eq:dwe} amounts to proving~\eqref{int:H}.
We therefore find, for~$t$ sufficiently large, \beno \Normt{\pd
  Uy}_\Ld\le C_0 e^{-C t}.\eeno We next differentiate
equation~\eqref{tcoue0a} with respect to~$t$ and insert~\eqref{tcoue0b}
to obtain \beno \pdd{u}{t}-\f\eta\rho \frac{\partial^2 }{\partial
  y^2}\brk{\pd u t} -\f{G}{\rho \lambda}\pdd{u}{y}=0.  \eeno The
function~$u$ satisfies the same equation as~$U$ and thus has the same
convergence rate. Applying the Gronwall Lemma to~\eqref{eq:tmt}
therefore implies, for~$t$ sufficiently large,
\begin{align}\label{d}
  \Normtp{\tau-\overline{\tau}}_\Ld\le C_0 e^{-Ct}.
\end{align}
Integrating~\eqref{tcoue0b} over~$\Omega$, we have \beno \lambda\f{d}{d
  t}\bar{\tau}=0, \eeno so that~$\bar{\tau}(t)=\overline{\tau_0}$ for
all times. We thus have the convergence estimate~\eqref{dr5}.
\endproof

\section{Longtime behaviour for non-homogeneous boundary conditions in a
  simple case}
\label{sec:nhom}
In this section, we study the longtime behaviour of the
system~\eqref{eq:couette} supplied with \emph{non-homogeneous} boundary
conditions~$u(t,0)=0$ and~$u(t,1)=a$ (where~$a$ is a constant scalar
different from zero and chosen positive, without loss of
generality,~$a>0$).

We denote~$(\uinf,\tauinf,\finf)$ a stationary state to the
system~\eqref{eq:couette}. We only consider the simplified case
\begin{align} \label{ineq:finf} f_\infty > 0 \mbox{ everywhere}.
\end{align}

The only stationary state that satisfies~\eqref{ineq:finf} is made
explicit in subsection~\ref{ssec:ss}. In subsection~\ref{ssec:stab}, we
show convergence in the longtime to this stationary state for small
initial perturbations.  In subsection~\ref{ssec:edo}, we study the
longtime behaviour for initial data that satisfy~$f_0>0$ without any
smallness conditions, but only in a simplified case that reduces
system~\eqref{eq:couette} to a system of ordinary differential
equations.

We do not state any result for the convergence to stationary states when
fluidity vanishes on some part of~$\Omega$.

\subsection{Stationary state}
\label{ssec:ss}
The following lemma makes precise the stationary state that satisfies
the condition~\eqref{ineq:finf}.
\begin{lem}[Stationary state]
  \label{lem:ss}
  Supply system~\eqref{eq:couette} with non-homogeneous boundary
  conditions~$\uinf(0)=0$ and~$\uinf(1)=a>0$.  The unique stationary
  solution~$(\uinf,\tauinf,\finf)$ in~$\brk{H^1(\Omega)}^3$
  satisfying~\eqref{ineq:finf} reads
  \begin{align} \label{def:ssn0} (u_\infty, \tau_\infty, f_\infty)(y) =
    \brk{a y, \f{\sqrt{1+4\nu\xi G a}+1}{2\xi}, \f{\sqrt{1+4\nu\xi G
          a}-1}{2\nu}} .
  \end{align}
\end{lem}

\begin{rmk}
  It is easy to extend the above result to stationary solutions
  $(\uinf,\tauinf,\finf)$ in~$ H^1(\Om)\times L^\infty(\Om)\times
  L^\infty(\Om)$ that satisfy~$\finf \nequiv 0$. Introducing $
  \Omega_\infty =\BRK{y \in \Om,\finf(y)>0} $ and $ \beta_\infty =
  \mbox{meas}(\Omega_\infty), $ the set of such stationary solutions
  reads
  \begin{align*}
    \brk{\pd{\uinf}{y},\tauinf,\finf}(y)=
    \begin{cases}
      \brk{ \dfrac{a}{\beta_\infty},\ \tau_L,\ \dfrac{-1+\xi\tau_L}{\nu}
      }
      & \mbox{on } \Omega_\infty\\
      \brk{ 0,\ \eta \dfrac{a}{\beta_\infty} + \tau_L,\ 0} & \mbox{on }
      \Omega \backslash \Omega_\infty,
    \end{cases}
  \end{align*}
  with $\tau_L = \dfrac{1}{2\xi} \brk{1+\sqrt{1+4 \nu \xi G
      a/\beta_\infty}}$.
\end{rmk}

\proof The stationary states~$(\uinf,\tauinf,\finf) : \Om \rightarrow
\R$ of the system~\eqref{eq:couette} that satisfy~\eqref{ineq:finf} are
solutions of the following system
\begin{subequations}\label{eq:coust}
  \begin{empheq}[left=\empheqlbrace]{align}
    0 &= \eta \pdd{\uinf}{y} + \pd{\tauinf}y,  \label{eq:cousta}\\
    0 &= G \frac{\partial \uinf}{\partial y} -\finf \tauinf,
    \label{eq:coustb}\\
    \finf &= \dfrac{-1+\xi |\tauinf|}{\nu}. \label{eq:coustc}
  \end{empheq}
\end{subequations}
We now show that such a steady state is unique and explicitly identify
it. Since $\tauinf\in H^1(\Omega)$,~\eqref{eq:cousta} shows that $\uinf$
belongs to~$H^2(\Omega)$. We integrate~\eqref{eq:cousta}
and~\eqref{eq:coustb} over~$\Omega$ and obtain
\begin{align}
  K &= \eta \pd{\uinf}{y} + \tauinf,
  \label{eq:K}\end{align}
where~$K$ is a constant and, using the boundary conditions on~$\uinf$,
\begin{align}
  \into \finf\tauinf &= G a.
  \label{eq:intob}\end{align}
We combine~\eqref{eq:coustb} and~\eqref{eq:K} to obtain
\begin{align}\label{eq:K2}
  \brk{\f\eta G \finf+1} \tauinf = K
\end{align}
so that~$\tauinf$ has the constant sign of
$K$. Equation~\eqref{eq:intob} then implies that~$\tauinf$, thus~$K$
are positive.\\
We now claim that~$\tauinf$ is constant over~$\Omega$:
inserting~\eqref{eq:coustc} in~\eqref{eq:K2}, we obtain that~$\tauinf$
satisfies
\begin{align*}
  \dps \tauinf \brk{ 1+ \eta \dfrac{-1+\xi \tauinf}{G\nu} } = K.
\end{align*}
It is easy to see that this equation has a unique positive
solution~$\tauinf$.  It follows from~\eqref{eq:K} that~$\pd \uinf y$ is
constant throughout~$\Omega$ so that, using the boundary conditions,
$\uinf(y)=a y.$ We rewrite equation~\eqref{eq:coustb} as
\begin{align*}
  G a = \dfrac{-1+\xi \tauinf}{\nu} \tauinf,
\end{align*}
to find the value of
\begin{align*}
  \tauinf=\dfrac{1}{2\xi} (1+\sqrt{1+4\nu \xi G a}).
\end{align*}
The stationary state reads~$\brk{a y, \tauinf, \dfrac{-1
    +\xi\tauinf}\nu}$, that is~\eqref{def:ssn0}.
\endproof

\subsection{Longtime behaviour with smallness assumption}
\label{ssec:stab}
The following theorem states the convergence in the longtime to the
stationary state~\eqref{def:ssn0} for small initial perturbations.

\begin{thm} \sloppy\label{th:stab} Supply system~\eqref{eq:couette} with
  non-homogeneous boundary conditions~$u(t,0)=0$ and~$u(t,1)=a>0$.
  Consider the solution $(u,\tau,f)$
  the existence and uniqueness of which have been established in
  Theorem~\ref{existence} and the associated stationary
  state~$(\uinf,\tauinf,\finf)$ defined by~\eqref{def:ssn0}. There
  exists~$\veps>0$ (sufficiently small so that at least~$\tau_0$ and
  $f_0$ are positive 
  ), such that, if the initial data~$(u_0, \tau_0, f_0)$
  for~\eqref{eq:couette} satisfy
  \begin{align*}
    \Norm{u_0-\uinf}_{H^1}^2+\Norm{\tau_0-\tau_\infty}_\Li^2
    +\Norm{f_0-f_\infty}_\Li^2\le \veps^2,
  \end{align*}
  then the solution~$(u, \tau, f)$ of system~\eqref{eq:couette}
  converges, as~$t$ goes to infinity, to~$(\uinf,\tauinf,\finf)$
  in~$H^1(\Omega)\times L^\infty(\Omega) \times L^\infty(\Omega)$.

  More precisely, there exist a constant~$C$ independent from~$\veps$,
  time and initial data and a constant~$\kgp_\veps$ independent from
  time such that, for~$t$ sufficiently large,
  \begin{align}\label{ineq:convs}
    \Norm{u(t,\cdot)-\uinf(\cdot)}_{H^1}
    +\Norm{\tau(t,\cdot)-\tau_\infty(\cdot)}_\Li
    +\Norm{f(t,\cdot)-f_\infty(\cdot)}_\Li \leq \kgp_\veps
    e^{-(C-\epsilon)t}.
  \end{align}
  \fussy\end{thm}

\begin{rmk} \label{rmk:h1} It is indeed possible, under the same
  assumptions, to prove that both~$\tau$ and~$f$ converge to zero
  in~$H^1(\Omega)$ and not only in~$L^\infty(\Omega)$. The proof is more
  tedious. We omit it here for brevity and refer to~\cite{these-david}.
\end{rmk}

Before we get to the proof , we note that we will return to
system~\eqref{eq:couette}, and not~\eqref{tcoue} since of course
boundary conditions will play a crucial role throughout the section.  We
also rewrite system~\eqref{eq:couette} as
\begin{equation}
  \left\{
    \begin{aligned}
      \rho \pd{u}{t} &= \eta \frac{\partial^2 }{\partial y^2}((u-\uinf)+\uinf)+ \frac{\partial }{\partial y} ((\tau-\tauinf)+\tau_\infty),\\
      \lambda \pd{\tau}{t} &= G \frac{\partial }{\partial y}((u-\uinf)+\uinf) - ((f-\finf)+f_\infty) ((\tau-\tauinf)+\tau_\infty), \\
      \pd{f}{t} &= (-1+\xi |(\tau-\tauinf)+\tau_\infty|)
      ((f-\finf)+f_\infty)^2-\nu ((f-\finf)+f_\infty)^3.
    \end{aligned}
  \right.\nonumber
\end{equation}
To lighten the notation, we henceforth denote~$(u,\tau,f)$ instead
of~$(u-\uinf,\tau-\tauinf,f-\finf)$ and consider
\begin{subequations}
  \begin{empheq}[left=\empheqlbrace]{align}
    \rho \pd{u}{t} &= \eta \pdd{u}{y}+  \pd\tau y, \label{pcouea}\\
    \lambda \pd{\tau}{t} & = G \pd{u}{y} -f_\infty\tau-\tau_\infty
    f-f\tau, \label{pcoueb} \\
    \pd{f}{t} & = -\nu(f+f_\infty)^2f+\xi(f+2f_\infty)f\tau+\xi
    f_\infty^2\tau, \label{pcouec}
  \end{empheq}
  \label{pcoue}
\end{subequations}
supplied with \emph{homogeneous} boundary conditions on~$u$ and initial
data that satisfy
\begin{align*}
  \Norm{u_0}_{H^1}^2+\Norm{\tau_0}_\Li^2+\Norm{f_0}_\Li^2\le \veps^2.
\end{align*}

\proof The proof is divided into three steps. The first step establishes
a priori estimates on system~\eqref{pcoue}.
In the second step, we show that the solution remains small for
sufficiently small perturbations. In Step~\ref{stabastab}, we show that,
still for small perturbations, the solution converges to the steady
state and that the rate of convergence is exponential.

As in the previous proofs,~$C$ and~$\kgp_\veps$ denote various constants
the value of which may vary from one instance to another, the actual
value being irrelevant.

\newstep{A priori energy estimates}{stab:ee} We argue as in
Step~\ref{s:fee} of the proof of Theorem~\ref{existence}.  We
multiply~\eqref{pcouea},~\eqref{pcoueb} and \eqref{pcouec} respectively
by~$u$,~$\tau$ and~$f$, integrate over~$\Omega$ and find
\begin{align*}
  \f{\rho}2\f{d}{d t}\Normt{u}_\Ld^2+\eta\Normt{\pd{u}{y}}_\Ld^2&=
  \into u \pd{\tau}{y}(t,\cdot),\\
  \f{\lambda}2\f{d}{dt}\Normt{\tau}_\Ld^2+\into
  (f+f_\infty)\tau^2(t,\cdot)&=
  G \into \tau \pd{u}{y}(t,\cdot)-\tau_\infty \into f\tau(t,\cdot),\\
  \f{1}2\f{d}{d t}\Normt{f}_\Ld^2+\nu \into (f+f_\infty)^2
  f^2(t,\cdot)&= \xi f_\infty^2 \into f\tau(t,\cdot)+\xi \into
  (f+2f_\infty)f^2\tau(t,\cdot).
\end{align*}
Combining these estimates leads to
\begin{align}\label{s1} \f12\f{d}{d t}&\brk{ \rho
    G\Normt{u}_\Ld^2+\lambda\Normt{\tau}_\Ld^2+\f{\tau_\infty}{\xi\finf^2}\Normt{f}_\Ld^2}+\eta
  G\Normt{\pd{u}{y}}_\Ld^2+f_\infty\Normt{\tau}_\Ld^2 \nonumber\\&
  +\f{\tau_\infty\nu}{\xi}\Normt{f}_\Ld^2 - C
  \Norm{f}_{L^\infty_T(L^\infty)}(\Norm{f
  }_{L^\infty_T(L^\infty)}+2\finf)\brk{\Normt{\tau}_\Ld^2+\Normt{f}_\Ld^2}
  \le 0 , \end{align} using the~$L^\infty$-estimate on~$[0,T]$ on~$f$
established in the proof of Theorem~\ref{existence} and the Young
inequality. We now use Step~\ref{s:apf} of the proof of
Theorem~\ref{existence} and more precisely the estimate~\eqref{en3}
on~$U$ defined by~\eqref{def:U}.  We have
\begin{align}\label{s2}
  \f{d}{d t}
  \Normt{\pd{U}{y}}_\Ld^2+\f{\eta}{\rho}\Normt{\pdd{U}{y}}_\Ld^2 \le
  C\brk{\Normt{u}_\Ld^2+
    \Norm{f}_{L^\infty_T(L^\infty)}^2\Normt{\tau}_\Ld^2}.
\end{align}

Now that we have estimates in Sobolev spaces, we turn to point wise
estimates on~$\tau$ and~$f$. We refine our argument in Step~\ref{s:li}
of the proof of Theorem~\ref{existence}. We rewrite the evolution
equation~\eqref{pcoueb} as
\begin{align*}
  \lambda \pd{\tau}{t} + \brk{f+\finf+\f{G}{\eta}}\tau=
  G\pd{U}{y}+\f{G}{\eta}\bar{\tau}-\tauinf f,
\end{align*}
multiply it by~$\tau$, apply the Young inequality and obtain
\begin{align}\label{eqp:tau}
  \f{\lambda}2\f{d}{d t}|\tau|^2+\brk{f+\finf+\f{G}{2\eta}}|\tau|^2\le
  \eta \Normt{\pd{U}{y}}_\Li^2+ \Normt{\tau}_\Ld^2 -\tauinf f \tau.
\end{align}
Similarly, we multiply~\eqref{pcouec} by~$f$ and find
\begin{align}\label{eqp:f}
  \f12\f{d}{d t}|f|^2+\nu (f+\finf)^2|f|^2=
  \xi(f+2f_\infty)\tau|f|^2+\xi f_\infty^2 f\tau.
\end{align}
We combine~\eqref{eqp:tau} and~\eqref{eqp:f} and use the Poincar\'e
inequality on~$\pd Uy$, the spatial average of which is zero, to obtain
\begin{align} \label{eq:fetau} \f{\lambda}2\f{d}{d t}&|\tau|^2 +
  \f{\tauinf}{2\xi\finf^2}\f{d}{d t}|f|^2
  +\brk{f+\finf+\f{G}{2\eta}}|\tau|^2 \nonumber\\&
  +\f{\tauinf}{\xi\finf^2} \brk{\nu(f+\finf)^2-
    \xi(\Norm{f}_{L^\infty_T(L^\infty)}+2f_\infty)\Norm{\tau}_{L^\infty_T(L^\infty)}}|f|^2
  \nonumber\\&\le \eta \Normt{\pdd{U}{y}}_\Ld^2+ \Normt{\tau}_\Ld^2 .
\end{align}

\newstep{Smallness of the solution for small perturbations}{stabstab} We
now prove that, for~$\epsilon \in (0,1)$ to be fixed later on and an
initial condition satisfying
\begin{align} \label{eq:petit0}
  \Norm{u_0}_{H^1}^2+\Norm{\tau_0}_\Li^2+\Norm{f_0}_\Li^2\le \epsilon^2,
\end{align}
we have, for all time~$t>0$,
\begin{align} \label{ineq:stab}
  \Normt{u}_{H^1}^2+\Normt{\tau}_\Li^2+\Normt{f}_\Li^2 \le \epsilon.
\end{align}
We argue by contradiction and suppose
\begin{align*}
  T_M= \inf\BRK{t\in\R^+|
    \brk{\Norm{u(t,\cdot)}_{H^1}^2+\Norm{\tau(t,\cdot)}_\Li^2+
      \Norm{f(t,\cdot)}_\Li^2}\ge \epsilon } \mbox{ is finite.}
\end{align*}
For all~$t\le T_M$, we use the estimates from the previous step.
For~$\epsilon$ sufficiently small such that all the terms in the
left-hand side of~\eqref{s1} are positive (this gives one condition
on~$\epsilon$), we have, integrating~\eqref{s1} from~$0$ to~$t$,
\begin{align*}
  \rho& G\Normt{u}_\Ld^2+\lambda\Normt{\tau}_\Ld^2+
  \f{\tau_\infty}{\xi\finf^2}\Normt{f}_\Ld^2 \nonumber\\&+\int_0^t \brk{
    \eta G\Norms{\pd{u}{y}}_\Ld^2+\frac\finf2 \Norms{\tau}_\Ld^2
    +\f{\tau_\infty\nu}{\xi}\Norms{f}_\Ld^2} ds\le C \epsilon^2.
\end{align*}
Integrating~\eqref{s2} from~$0$ to~$t$ then yields
\begin{align} \label{ineqf:U}
  \Normt{\pd{U}{y}}_\Ld^2+\int_0^t\Norms{\pdd{U}{y}}_\Ld^2ds \le
  C\epsilon^2 .
\end{align}
We now integrate~\eqref{eq:fetau} from~$0$ to~$t$ and get
\begin{align}\label{ineqf:fetau}
  \f{\lambda}2|\tau|^2 +\f{\tauinf}{2\xi\finf^2} |f|^2
  \le C\epsilon^2.
\end{align}
For all~$t\le T_M$,~\eqref{ineqf:U} and~\eqref{ineqf:fetau} imply
\begin{align*}
  \Normt{\pd{u}{y}}_\Ld^2+ \Normt{\tau}_\Li^2 + \Normt{f}_\Li^2 \le
  C\epsilon^2.
\end{align*}
Choosing~$\epsilon$ sufficiently small such that $C\epsilon^2<\epsilon$
(which gives another condition on~$\epsilon$) contradicts the definition
of~$T_M$, and so~$T_M=\infty$. It follows that~\eqref{ineq:stab} holds
for all time~$t>0$, the solution remains small.

\newstep{Convergence to the stationary state}{stabastab} We now prove
that, if the initial data satisfy~\eqref{eq:petit0}, then the solution
converges exponentially fast to the stationary state in the longtime.
For~$t$ sufficiently large,~\eqref{s1} implies that
\begin{align} \label{eq:assl2}
  \Normt{u}_\Ld^2+\Normt{\tau}_\Ld^2+\Normt{f}_\Ld^2\le \kgp_\veps
  e^{-(C-\epsilon)t}.
\end{align}
Adding~\eqref{s2} multiplied by~$2\rho$ to~\eqref{eq:fetau} leads to
\begin{align*}
  \f{d}{d t}&\brk{2\rho\Normt{\pd{U}{y}}_\Ld^2+\f{\lambda}2|\tau|^2 +
    \f{\tauinf}{2\xi\finf^2}|f|^2} +\eta\Normt{\pdd{U}{y}}_\Ld^2
  +(f+\finf+\f{G}{2\eta})|\tau|^2 \nonumber\\& +\f{\tauinf}{\xi\finf^2}
  \brk{\nu(f+\finf)^2-
    \xi(\Norm{f}_{L^\infty_T(L^\infty)}+2f_\infty)\Norm{\tau}_{L^\infty_T(L^\infty)}}|f|^2
  \nonumber\\&\le C\brk{\Normt{u}_\Ld^2+ \Normt{\tau}_\Ld^2
    +\Normt{f}_\Ld^2+\Norm{f}_{L^\infty_T(L^\infty)}^2\Normt{\tau}_\Ld^2}.
\end{align*}
We use the Poincar\'e inequality on~$\pd Uy$, the spatial average of
which is zero, apply the Gronwall Lemma, insert~\eqref{eq:assl2} and
find
\begin{align*}
  \Normt{\pd{U}{y}}_\Ld^2+\Normt{\tau}_\Li^2+\Normt{f}_\Li^2 \le
  \kgp_\veps e^{-(C-\epsilon)t}.
\end{align*}
This convergence estimate is equivalent to~\eqref{ineq:convs} and we
have exponential convergence.

\endproof
\setcounter{step}{0}

\subsection{Longtime behaviour without smallness assumption (simplified
  case)}
\label{ssec:edo}
We now examine the longtime behaviour of system~\eqref{eq:couette}
supplied with \emph{not necessarily small} initial data ~$(u_0, \tau_0,
f_0)$.  We are unable to prove a general result and focus our attention
to the particular case where the initial condition is~$u_0=a y$ ($a$
positive constant), $\tau_0=\mbox{constant}=\overline{\tau_0}$,
$f_0=\mbox{constant}=\overline{f_0}>0$. In such a case, a substantial
simplification occurs. Indeed,~\eqref{eq:couette} reduces to the
following system of ordinary differential equations:
\begin{subequations}
  \begin{empheq}[left=\empheqlbrace]{align}
    \lambda \pd{\tau}{t} &= - f \tau+Ga \label{eq:odea}\\
    \pd{f}{t} &= (-1+\xi |\tau|) f^2-\nu f^3,\label{eq:odeb}
  \end{empheq}
  \label{eq:ode}
\end{subequations}
supplied with initial conditions~$\tau_0, f_0 \in\R$ with~$f_0>0$.\\
System~\eqref{eq:ode} has a unique steady state such that~$\finf > 0$
and it reads
\begin{align} \label{def:ssode} (\tauinf, \finf)=\brk{\f{\sqrt{1+4\nu\xi
        G a}+1}{2\xi}, \f{\sqrt{1+4\nu\xi G a}-1}{2\nu}}.
\end{align}
Indeed, such a steady state~$(\tauinf, \finf)$ satisfies~$\finf \tauinf=
G a$ (so that~$\tauinf>0$) and~$\nu\finf={-1+\xi\tauinf}$. Combining
these equations implies~$\brk{-1+\xi\tauinf} \tauinf = \nu G a$. This
equation has a unique solution given in~\eqref{def:ssode}.

In addition, we introduce the condition
\begin{align}\label{as}
  -\f1\lambda-2+2\xi(1+Ga)\brk{
    \f1{\sigma}+\f{\lambda\xi}{2Ga}\brk{\f{\nu\sigma+1}{\xi}+\f4{\xi}}^2
  }<0,
\end{align}
with
\begin{align}\label{def:sigma}
  \sigma=\min\BRK{\f{3Ga}{Ga\nu+4\tauinf}, \f{\sqrt{1+4\nu\xi G
        a}-1}{3\nu}}.
\end{align}
We are unable to perform our proof without this additional
assumption. The numerical simulations we perform (see
Figure~\ref{fig:_g1_dt5em3_raff}) however show convergence holds even
when~\eqref{as} is not satisfied.
 
\begin{thm}\label{th:ode}
  Supply system~\eqref{eq:ode} with initial conditions~$\tau_0, f_0
  \in\R$
  with~$f_0>0$.  Then the solution~$(\tau, f)$  remains bounded.\\
  In addition, under assumption~\eqref{as}, the solution~$(\tau, f)$
  converges to~$(\tauinf,\finf)$ in the longtime and the rate of
  convergence is exponential: for~$t$ sufficiently large,
  \begin{align}
    |\tau(t)-\tauinf|+|f(t)-\finf|\le C_0 e^{-C_r t},
  \end{align}
  where~$C_0$ is a constant independent from time and~$C_r$ reads
  \begin{align} \label{def:Cl} C_r=\begin{cases} \frac12
      \brk{\frac\finf\lambda+\nu\finf^2}-\frac12
      \sqrt\Delta, & \mbox{if}\quad \Delta \geq 0,\\
      \frac12 \brk{\frac\finf\lambda+\nu\finf^2}, & \mbox{if}\quad
      \Delta< 0,
    \end{cases}
  \end{align}
  with
  \begin{align} \label{def:delta} \Delta=\finf^2 \brk{
      \brk{\frac1\lambda+\nu\finf}^2 -4 \brk{\frac\nu\lambda \finf
        +\frac\xi\lambda \tauinf }}.
  \end{align}
\end{thm}

\proof The proof is divided into seven steps. Step~\ref{s:edos1}
introduces simplifications on the initial data and the system, that are
not restrictive for the longtime behaviour.  Some notation is given in
Step~\ref{s:edon}. A lower bound on~$f$ is derived in Step~\ref{s:edol}
and is used in Step~\ref{s:edob} to prove that the solution is
bounded. Further restrictions are made in Step~\ref{s:edos2} still
without loss of generality. Step~\ref{s:edoc} establishes the
convergence, which is proven to be exponential in Step~\ref{s:linstab}.

We consider until Step~\ref{s:edob} the maximal solution
to~\eqref{eq:ode} although the solution a posteriori exists for all
times because of boundedness.

\newstep{Simplifications on the initial data}{s:edos1} We show
that~$\tau$ and~$f$ solution to~\eqref{eq:ode} remain positive, possibly
after some time for~$\tau$.  We first remark that, since~$f_0>0$,~$f>0$
for all times, arguing as in Step~\ref{s:pos} of the proof of
Theorem~\ref{existence}.  On the other hand, if~$\tau\le0$ on some time
interval, evolution equation~\eqref{eq:odea} thus implies that~$\tau$
increases strictly on this time interval (recalling that~$a>0$). Hence,
there exist a time~$T_0$ such that $\tau(T_0)>0$. Moreover, for
all~$t>T_0$,~$\tau$ remains positive (since if~$\tau$ is zero at one
time~$T_1>T_0$, $\deriv\tau t(T_1)=Ga>0$, which is in contradiction
with~$\tau>0$ for~$t<T_1$.)

For the purpose of studying the longtime limit, we may always consider,
without loss of generality, the system
\begin{subequations}
  \begin{empheq}[left=\empheqlbrace]{align}
    \lambda \pd{\tau}{t} &= - f \tau+Ga \label{eq:podea}\\
    \pd{f}{t} &= (-1+\xi \tau) f^2-\nu f^3,\label{eq:podeb}
  \end{empheq}
  \label{eq:pode}
\end{subequations}
supplied with positive initial conditions~$\tau_0, f_0$.

\newstep{Some notation}{s:edon}

\begin{figure} [ht]
  \centering \scalebox{0.5}{\input{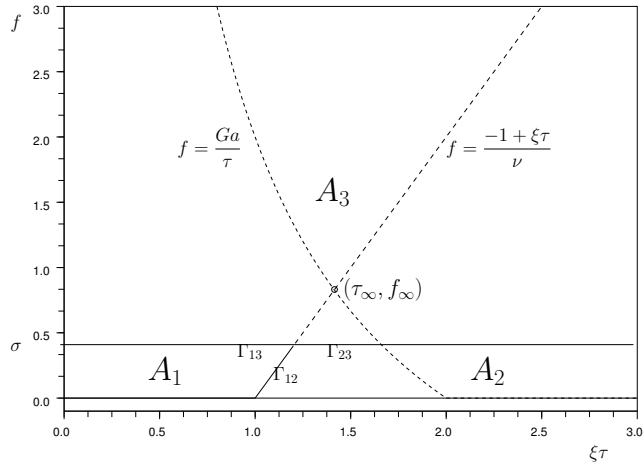}}
  \caption{Notation on~$(0,+\infty) \times (0,+\infty)$}
  \label{fig:edo_regions}
\end{figure}

{\sloppy We consider the three subdomains:
  \begin{align*}
    A_1&= \BRK{(\tau,f)|f\ge\frac{\xi \tau-1}{\nu},~f\le \sigma},\\
    A_2&=\BRK{( \tau, f)|f\le\frac{\xi \tau-1}{\nu},~f\le \sigma},\\
    A_3&= \BRK{(\tau,f)|f\ge \sigma} ,
  \end{align*}
}\fussy where we recall that~$\sigma$ is defined by~\eqref{def:sigma}.
We also introduce their intersections:
\begin{align*}
  \Gamma_{13} &= \BRK{(\tau, \sigma)|\tau\le \frac{\xi\sigma-1}{\nu}},\\
  \Gamma_{12}&= \BRK{(\tau,f)|f=\frac{\xi \tau-1}{\nu}, f\le \sigma},\\
  \Gamma_{23}&=\BRK{( \tau, \sigma)|\tau\ge\frac{\xi\sigma-1}{\nu}}.
\end{align*}
See Figure~\ref{fig:edo_regions} for a graphical description.

\newstep{Lower bound on~$f$}{s:edol} We now establish a lower bound for
the fluidity~$f$ in each domain.  In the cases~$(\tau_0, f_0)\in A_2$
or~$(\tau_0, f_0)\in A_3$, we have
\begin{align*}
  f\ge \min\BRK{f_0, \sigma }.
\end{align*}
The case~$(\tau_0, f_0)\in A_1$ requires more developments.  The
evolution equations~\eqref{eq:podea} and~\eqref{eq:podeb} respectively
rewrite
\begin{align*}
  \frac{d}{d t}\f1f &= 1-\xi \tau +\nu f,\\
  \f\lambda2 \frac{d}{d t} \brk{\tau-\f4\xi}^2 &=
  -f\tau\brk{\tau-\f4\xi} + Ga \brk{\tau-\f4\xi}.
\end{align*}
We combine these two equations and obtain
\begin{align}
  \lambda\xi \frac{d }{d t}
  \brk{\f{Ga}{\lambda\xi}\f1{f}+\f12\brk{\tau-\f4{\xi}}^2}&=-{3Ga}+
  {Ga\nu}f
  -\xi f\tau^2 +4f\tau\nonumber\\
  &\le -{3Ga}+\brk{{Ga\nu}+4\tauinf}\sigma\nonumber\\
  &\le 0 ,\label{ineq:subtil}
\end{align}
where we have used firstly that~$0\le f\le\sigma$ and~$\tau\le\tauinf$
in~$A_1$ and secondly~\eqref{def:sigma}.
Integrating~\eqref{ineq:subtil} yields
\begin{align*} \f{Ga}{\lambda\xi}\f1{f(t)}\le
  \f{Ga}{\lambda\xi}\max\BRK{\f1{f_0},
    \f1{\sigma}}+\f12\brk{\f{\nu\sigma+1}{\xi}+\f4{\xi}}^2,
\end{align*}
using that~$\tau\le \dfrac{\nu\sigma+1}{\xi}$.  We introduce
\begin{align}\label{def:m_f}
  m_f=\brk{\max\BRK{\f1{f_0}, \f1{\sigma}}+\f{\lambda\xi}{2Ga}
    \brk{\f{\nu\sigma+1}{\xi}+\f4{\xi}}^2 }^{-1},
\end{align}
which is therefore a lower bound for~$f$ in the region~$A_1$.  This
lower bound also holds for initial conditions that belong to~$A_2$
and~$A_3$ and we thus have, for all~$t>0$,
\begin{align}\label{lo} f(t)\ge m_f, \end{align}
with~$m_f$ defined by~\eqref{def:m_f}.

\newstep{Boundedness}{s:edob} The purpose of this step is to prove that
the solution~$(\tau,f)$
remains bounded.\\
Applying the Duhamel formula on~\eqref{eq:podea} yields
\begin{align*}
  \tau(t)=e^{ -\int_0^t \frac{f(s)}{\lambda}ds} \tau_0
  +\frac{Ga}{\lambda}\int_0^t e^{ -\int_s^t \frac{f(s')}{\lambda}ds'}
  ds,
\end{align*}
so that, using the lower bound~\eqref{lo} on~$f$,
\begin{align*}
  \tau(t)\le e^{-\dfrac{m_ft}\lambda} \tau_0+\f{Ga}{m_f}.
\end{align*}
Therefore,~$\tau$ is bounded, and there exists a time~$t_0$ such that,
for all~$t>t_0$,
\begin{align} \label{lo1} \tau(t)\le
  \f{Ga+1}{m_f}.
\end{align}
We now turn to the boundedness of~$f$. We
introduce~$M_\tau=\dfrac{Ga+1}{m_f}$ and~$M_f=\dfrac2\nu\brk{-1+\xi
  M_\tau}$.  We will show that, for all~$t>t_0$,
\begin{align}\label{ub:f}
  f(t) < \max\brk{f(t_0),M_f}.
\end{align}
We distinguish between two cases. Let us first suppose that~$\pd
ft(t_0)\ge0$.  In this case, $f(t_0)\le \dfrac1\nu
(-1+\xi\tau(t_0))<M_f$ because of~\eqref{lo1}.  Moreover, for
all~$t>t_0$, $f(t)<M_f$.  Indeed, by contradiction, if
\begin{align*}
  t_1=\inf\BRK{t>t_0,f(t_1)=M_f}<+\infty,
\end{align*}
then, by continuity, $\pd ft(t_1)\ge0$. On the other hand, we have
\begin{align*}
  \pd{f}{t}(t_1) &= f^2(t_1) (-1+\xi\tau(t_1)-\nu f(t_1))\\
  &< f^2(t_1) (-1+\xi M_\tau -\nu M_f)\\
  &<0,
\end{align*}
hence the contradiction. \\
In the other case~$\pd ft(t_0)<0$, $f$ strictly decreases until
(possibly) equality occurs at a later time~$t_3\ \brk{\pd ft(t_3)=0}$,
which leads to the previous case with~$t_3$ instead of~$t_0$. In any
case, we have obtained~\eqref{ub:f} for all~$t>t_0$.

\newstep{Further simplifications on the initial data}{s:edos2} Table~1
first summarizes how~$(\tau, f)$ behaves when it touches an intersection
line.
\begin{table}[ht]
  \label{tab:line}
  \centering
  \begin{tabular}{|c|c|c|c|}
    \hline
    starting line  & ~$\pd f t$  & ~$\pd\tau t$  &  entering region  \\
    \hline
    $\Gamma_{13}$    &  -          &  +            &  $A_1$                \\
    $\Gamma_{12}$    &  0          &  +            &  $A_2$               \\
    $\Gamma_{23}$    &  +          &               &  $A_3$              \\
    \hline
  \end{tabular}
  \caption{Motion on intersection lines}
\end{table}
We use Table~1 to show that the solution enters region~$A_3$ at some
time.

In region~$A_1$, we have~$\pd{\tau}{t}>0$, so that there does not exist
any periodic orbit inside region~$A_1$. There is also no steady state in
this region. Using the Poincar\'e-Bendixson Theorem on the bounded
solution of ordinary differential equation system~\eqref{eq:pode}, the
solution leaves
region~$A_1$ at some time. According to Table~1, it enters region~$A_2$.\\
Applying similar arguments on region~$A_2$ where~$\pd{f}{t}>0$, the
solution enters region~$A_3$ at some time.\\
We can therefore restrict the studying of the longtime limit to initial
data~$(\tau_0, f_0)$ that belongs to region~$A_3$, without loss of
generality.

The bounds~\eqref{lo} and~\eqref{lo1} become
\begin{align}\label{lo2} f(t)\ge \brk{ \f1{\sigma}+
    \f{\lambda\xi}{2Ga}\brk{\f{\nu\sigma+1}{\xi}+\f4{\xi}}^2
  }^{-1}.  \end{align} and \begin{align}\label{lo3} \tau(t)\le
  (Ga+1)\brk{ \f1{\sigma}+\f{\lambda\xi}{2Ga}
    \brk{\f{\nu\sigma+1}{\xi}+\f4{\xi}}^2 }.
\end{align}

\newstep{Convergence}{s:edoc} We introduce
\begin{align*}
  G(\tau, f)=-\f1{\lambda} f \tau+\f1{\lambda}Ga
\end{align*}
and
\begin{align*}
  F(\tau,f) = (-1+\xi \tau) f^2-\nu f^3.
\end{align*}
We have
\begin{align*}
  \frac{\pa G(\tau, f)}{\pa \tau}+\frac{\pa F(\tau, f) }{\pa
    f}&=f\brk{-\f1\lambda-2+2\xi \tau-3\nu f},\\
  &<f\brk{ -\f1\lambda-2+2\xi(1+Ga)\brk{
      \f1{\sigma}+\f{\lambda\xi}{2Ga}\brk{\f{\nu\sigma+1}{\xi}+\f4{\xi}}^2
    }},
\end{align*}
using~\eqref{lo3} and the positivity of~$f$. Because of our
assumption~\eqref{as}, the right-hand sides~$F$ and~$G$
of~\eqref{eq:pode} satisfy
\begin{align*}
  \frac{\pa G(\tau, f)}{\pa \tau}+\frac{\pa F(\tau, f) }{\pa f}<0.
\end{align*}
According to the Dulac Criterion, there does not exist any periodic
orbit for~\eqref{eq:pode}. Since it has only one steady state~$(\tauinf,
\finf)$, the solution converges to it:
\begin{align*}
  \lim_{t\rightarrow\infty}\brk{|\tau(t)-\tauinf|+|f(t)-\finf|}=0.
\end{align*}

\newstep{Exponential convergence}{s:linstab} Now that we have
convergence to the steady state, we can use linear
stability. System~\eqref{eq:pode} linearized around the stationary
state~$(\tauinf,\finf)$ reads
\begin{align*} \frac{d }{d t} \brk{\begin{array}{cc} \tau_l \\
      f_l \end{array}}=\brk{\begin{array}{cc}
      -\frac\finf\lambda&-\frac\tauinf\lambda \\ \xi \finf^2& -\nu
      \finf^2\end{array}} \brk{\begin{array}{cc} \tau_l \\
      f_l \end{array}}.
\end{align*}
The eigenvalues of the associated matrix depend on the sign of~$\Delta$
defined by~\eqref{def:delta}.  If~$\Delta<0$, the eigenvalues are
complex and their real part is ~$-\frac12
\brk{\frac\finf\lambda+\nu\finf^2}$. If~$\Delta \geq 0$, the eigenvalues
are real negative, the smaller one in absolute value is~$-\frac12
\brk{\frac\finf\lambda+\nu\finf^2}+\frac12 \sqrt\Delta~$.  The real part
of the eigenvalues gives the rate of convergence and hence of values
of~$C_r$ in~\eqref{def:Cl}.

\setcounter{step}{0}
\endproof

\section{Numerical results}

In this section, we present numerical simulations that complement the
theoretical results on the behaviour of the previous sections.

We simulate numerically~\eqref{eq:couette} in the
interval~$\Omega=[0,1]$ and the interval~$[0,T]$ for~$T=10000$. The
system is supplied either with \emph{homogeneous} boundary conditions or
non-homogeneous boundary conditions~$u(t,0)=0$ and~$u(t,1)=a$ for all
time~$t\in[0,T]$. In the latter case, we take~$a=1$.  As for the initial
conditions, we take sinusoidal functions for all three fields. The
values of~$u_0$ oscillate between~$-0.002$ and~$0.002$ for
\emph{homogeneous} boundary conditions and between~$0$ and~$a$
otherwise.  The values of~$\tau_0$ and~$f_0$ oscillate between~$-0.5$
and~$0.5$.

We use the following set of physical parameters. The
density~$\rho=0.001$ and the viscosity~$\eta=1$ so that the Reynolds
number is low. The elastic modulus~$G$ and the coefficients~$\xi$
and~$\nu$ are equal to one. The characteristic relaxation time~$\lambda$
is~$0.5$ unless otherwise stated.
  
System~\eqref{eq:couette} is solved using a constant time step~$\Delta t
= 0.005$ with the following time scheme:
\begin{subequations}
  \begin{empheq}[left=\empheqlbrace]{align}
    \frac\rho{\Delta t} \brk{u_n-u_{n-1}} &= \eta \pdd{u_n}{y} +
    \pd{\tau_{n-1}}{y},\label{eq:couettena}\\
    \frac\lambda{\Delta t} \brk{\tau_n-\tau_{n-1}}  &= G \pd{u_n}{y} - f_{n-1} \tau_{n-1},\label{eq:couettenb} \\
    \frac1{\Delta t} \brk{f_n-f_{n-1}} &= (-1+\xi |\tau_n|) f_{n-1}
    f_n-\nu f_{n-1} f_n^2.\label{eq:couettenc}
  \end{empheq}
  \label{eq:couetten}
\end{subequations}
For the space variable, we use linear $\mathbb{P}1$ finite elements
for~$u$ and piecewise constant finite elements for both~$\tau$
and~$f$. Note that, in contrast to the approximating
system~\eqref{eq:atcoue} we used for our theoretical proof, we
take~$\tau_{n-1}$ instead of~$\tau_n$ in the right-hand sides
of~\eqref{eq:couettena} and~\eqref{eq:couettenb}. This allows us to
solve each equation separately. This choice is made for
simplicity. Other approaches could have been employed. For our tests, we
use elements of constant size~$h=0.002$ and perform the computations
using Scilab\cite{Scilab}.

\subsection{\emph{Homogeneous} boundary conditions}
\label{ssec:numh}
\begin{figure} [ht]
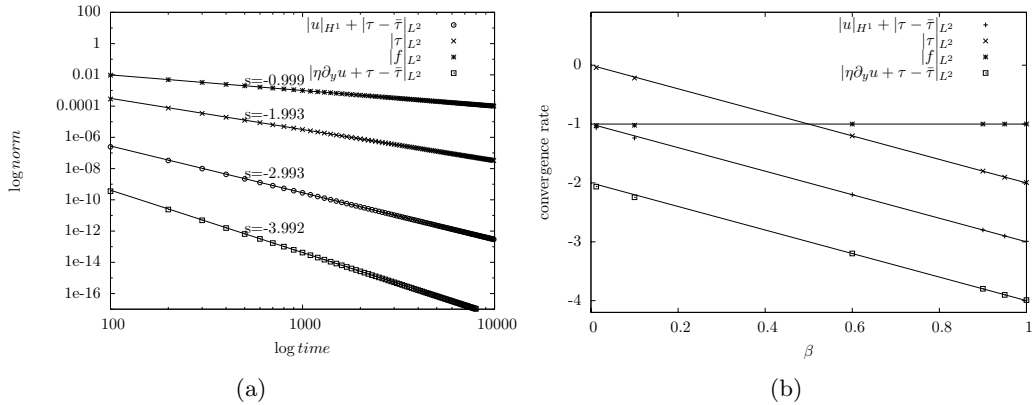

  \centering \subfigure[
  ]{\label{fig:_bc0}\scalebox{0.55}{\input{\inc/__bc0.tex}}} \subfigure[
  ]{\label{fig:beta}\scalebox{0.55}{\input{\inc/_beta.tex}}}
  \caption{(a) Time evolution in log-log scale for \emph{homogeneous}
    boundary conditions; the points are the simulated trajectories; the
    lines and the corresponding slopes~$s$ are fitted. (b) Fitted
    convergence rates~$s$ for~$\beta=0, 0.01, 0.1, 0.6, 0.9, 0.99$; the
    lines are the theoretical convergence rates function of~$\beta$.}
\end{figure}

We first focus on the \emph{homogeneous} boundary conditions on~$u$
considered in Section~\ref{sec:hom}.  The case~$f_0\equiv 0$, that
implies~$f\equiv 0$ for all times, is uninteresting numerically. We
therefore only show results for~$f_0 \nequiv 0$. In this case, we have
convergence to the stationary state~$(0,0,0)$ as proven in
Theorem~\ref{th:longhom1}. The convergence estimates are established in
Theorem~\ref{th:longhomt}.  We recall the parameter
\begin{align*}
  \beta=meas\BRK{y\in \Omega |f_0(y)> 0}
\end{align*}
and these convergence rates: for~$\vveps$ arbitrarily small, there exist
various constants~$\kgp_\vveps$ and a time~$t_0$ , such that, for
all~$t>t_0$,
\begin{align*}
  \Normt{\tau}_\Ld&\le
  \kgp_\vveps (1+t)^{-\f{\beta}{\lambda}(1-\vveps)},\\
  \f1{\f1{f(t_0,y)}+(1+\vveps)(t-t_0)}\le f(t,y)&\le
  \f{1}{\f1{f(t_0,y)}+
    (1-\vveps)(t-t_0) }\\
  \Normt{u}_{H^1}+\Normtp{\tau-\overline{\tau}}_\Ld&\le
  \kgp_\vveps (1+t)^{-1-\f{\beta}{\lambda}(1-\vveps)}\\
  \Normtp{\eta \pd{u}{y}+\tau-\overline{\tau}}_\Ld&\le \kgp_\vveps
  (1+t)^{-2-\f{\beta}{\lambda}(1-\vveps)}.
\end{align*}
Note that the last three estimates are exactly the same as in
Theorem~\ref{th:longhomt}, the first estimate is an immediate
consequence of~\eqref{dr2} and~\eqref{dr3}.  We now check that these
estimates are sharp.  We begin with the case~$f_0>0$ on~$\Omega$ that
is~$\beta=1$. The evolutions of $ \Normt{\tau}_\Ld, \Normt{f}_\Ld,
\Normt{u}_{H^1}+\Normtp{\tau-\overline{\tau}}_\Ld, \Normtp{\eta
  \pd{u}{y}+\tau-\overline{\tau}}_\Ld $ are represented in
Figure~\ref{fig:_bc0}. We use a log-log representation.  The slopes~$s$,
which correspond to a decrease as~$t^s$, are fitted on the numerical
results and indicated on Figure~\ref{fig:_bc0}: the numerical
convergence rates, obtained with~$\lambda=0.5$, are in good agreement with the estimates.\\
We next consider cases where~$f_0=0$ on some part of the domain. In
Figure~\ref{fig:beta}, we show simulations obtained with different
values of~$\beta$. For each simulation, that is for each value
of~$\beta$ considered, the convergence rates are fitted and represented
as a function of~$\beta$. The
numerical and theoretical convergence rates~$s$ agree. \\
We have extended these results to the other values of~$\lambda$ than
$\lambda=0.5$ and other values of the parameters $\rho, \eta, G, \xi,
\nu$ to check that the convergence estimates of
Theorem~\ref{th:longhomt} depend only on~$\lambda$ and~$\beta$ and are
indeed sharp.

\subsection{Non-homogeneous boundary conditions}
\begin{figure} [ht]
  \centering \subfigure
  {\scalebox{0.5}{\input{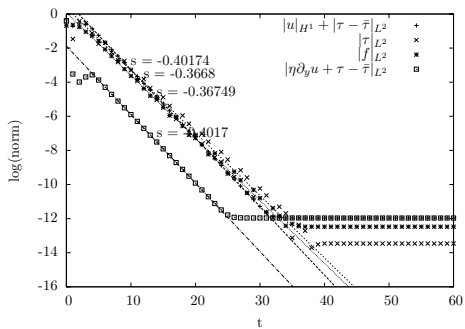}}}
  \caption{Time evolution of the perturbation to equilibrium in
    semi-logarithmic scale for non-homogeneous boundary conditions ; the
    points are the simulated trajectories and the line and the
    corresponding slope are fitted.}
  \label{fig:_g1_dt5em3_raff}
\end{figure}

The longtime behaviour for non-homogeneous boundary conditions has been
studied in Section~\ref{sec:nhom}. We consider only stationary
states~$(\uinf,\tauinf,\finf)$ that satisfy~$\finf>0$ everywhere. We
have shown that such a steady state~\eqref{def:ssn0} is unique. We
established in Theorem~\ref{th:stab} that we have convergence to this
steady state for small perturbations.  To have convergence, we of course
need to assume~$f_0> 0$ everywhere. We observe numerically that no other
condition, and specifically non assumption on the smallness of the data,
is required.  We consider the
perturbations~$(u-\uinf,\tau-\tauinf,f-\finf)$ to equilibrium and show
that they vanish in the longtime, see Figure~\ref{fig:_g1_dt5em3_raff}.
The evolution is plotted in semi-logarithmic scale.  The convergences of
the various norms $ \Normt{\tau}_\Ld, \Normt{f}_\Ld,
\Normt{u}_{H^1}+\Normtp{\tau-\overline{\tau}}_\Ld$, $\Normtp{\eta
  \pd{u}{y}+\tau-\overline{\tau}}_\Ld$ are indeed exponential.

\begin{figure} [ht]
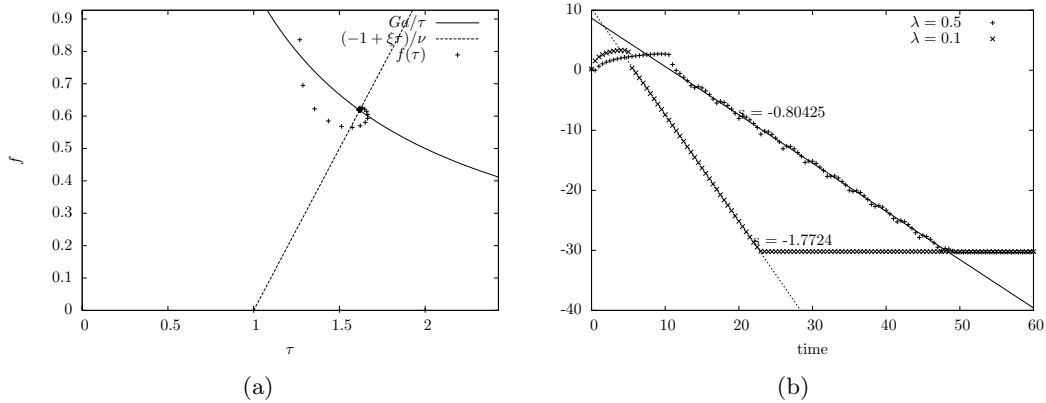

  \centering \subfigure[
  ]{\label{fig:Dstau12s}\scalebox{0.55}{\input{\inc/_Dstau12s.tex}}}
  \subfigure[
  ]{\label{fig:ode_decay}\scalebox{0.55}{\input{\inc/_ode_decay.tex}}}
  \caption{For the system of ordinary differential
    equations~\eqref{eq:ode}, time evolution (a) in the space~$(\tau,f)$
    for~$\lambda=0.5$; (b) of the perturbation to equilibrium in
    semi-logarithmic scale for~$\lambda=0.5,0.1$; the points are the
    simulated trajectories and the line and the corresponding slope are
    fitted. }
\end{figure}

In section~\ref{sec:nhom}, in order to establish a result without any
smallness assumption, we have considered a particular initial data that
reduces~\eqref{eq:couette} to the ordinary differential equation
system~\eqref{eq:ode}. We have obtained convergence to the stationary
state~\eqref{def:ssode} and explicit formula for the rate of
convergence.  Numerically, we observe convergence even when the
condition~\eqref{as}, which was assumed for the proof of
Theorem~\ref{th:ode}, is not satisfied.  The time evolution is shown in
the space of~$(\tau,f)$ and the convergence are represented in
Figure~\ref{fig:Dstau12s}.  We check that the convergence is exponential
as observed numerically in the general case of~\eqref{eq:couette} (see
Figure~\ref{fig:_g1_dt5em3_raff}). Moreover, we compute the convergence
rate and compare it to the theoretical rate~$C_r$ defined
by~\eqref{def:Cl}.  The evolution of the perturbation
function~$|\tau(t)-\tauinf|+|f(t)-\finf|$ is plotted as a function of
time in semi-logarithmic scale in Figure~\ref{fig:ode_decay}.  The first
case~$\lambda=0.5$ correspond to the case when the eigenvalues of the
associated linearized system are complex, the expected value of~$C_r$
is~$0.8090$; the other case~$\lambda=0.1$ is when the eigenvalues are
real negative, the expected value of~$C_r$ is~$1.7895$. The theoretical
and numerical value agree.

\acknowledgements{The authors are grateful to Fran\c cois Lequeux (ESPCI
  Paris) for many stimulating and enlightening discussions on the
  rheology of complex fluids.}

\end{document}